\documentclass[reqno]{amsart}
\usepackage{amsmath, amssymb, amscd}

\numberwithin{equation}{section}
\theoremstyle{plain}
\newtheorem{thm}{Theorem}[section]
\newtheorem{cor}[thm]{Corollary}

\newtheorem{lem}[thm]{Lemma}
\newtheorem{prop}[thm]{Proposition}
\newtheorem{Def}[thm]{Definition}
\theoremstyle{remark}
\newtheorem{rem}[thm]{Remark}

\title{Koornwinder polynomials and affine Hecke algebras}
\author{Jasper V. Stokman}

\address{Jasper V. Stokman, 
KdV Institute for Mathematics, Universiteit van Amsterdam,
Plantage Muidergracht 24, 1018 TV Amsterdam, The Netherlands.}
\email{jstokman@science.uva.nl}
\date{April 18th, 2005\\
\indent 2000 {\it Mathematics Subject Classification}. 33D52, 33D80.}
\begin{document}

\keywords{(Non-)symmetric Koornwinder polynomial, affine root system,
(double) affine Hecke algebra, Noumi representation,
Cherednik-Dunkl type difference-reflection operator,
intertwiner, (bi-)orthogonality, quadratic norm,
evaluation formula, duality}


\begin{abstract}
In this paper we derive the 
bi-orthogonality relations, diagonal term evaluations
and evaluation formulas for the non-symmetric
Koornwinder polynomials. 
For the derivation we use 
certain representations of the (double) affine Hecke algebra 
which were originally defined by Noumi and Sahi. 
The structure of the diagonal terms is clarified
by expressing them as residues of the 
bi-orthogonality weight function.
We  furthermore give the explicit connection between 
the non-symmetric and the (anti-)symmetric theory. 
\end{abstract}

\maketitle



\section{Introduction}

Cherednik \cite{C1}--\cite{C5} and Macdonald \cite{M3}
clarified the structure of Macdonald polynomials
using certain representations of affine Hecke algebras in terms
of difference-reflection operators.
The underlying data stem from a fixed, reduced, irreducible
root system $\Sigma$. The
degrees of freedom are, besides the deformation parameter $q$,
the number of different root length occurring in $\Sigma$ (so at most two).
In fact, their work shows that the Macdonald polynomials are
naturally attached to the reduced,
affine root system $\widetilde{\Sigma}$ associated with $\Sigma$. 

As announced by Macdonald 
\cite[sect. 8]{C3}, and partially carried out by Noumi
\cite{N}, Sahi \cite{Sa}, \cite{Sa2}, Noumi \& Stokman \cite{NS} 
and Nishino et al \cite{NUKW},
the theory naturally extends to the setting of arbitrary
(not necessarily reduced) irreducible affine root systems. This in particular
allows to incorporate the very general six parameter family of Koornwinder
\cite{Ko} polynomials in the theory.

The affine root system underlying the 
Koornwinder polynomials is the non-reduced, irreducible affine
root system $S$ of type $C^\vee C_n$, which was introduced by Macdonald in
\cite{M1}. The affine root
system $S$ contains all (possibly non-reduced) irreducible
affine root systems of classical type as an affine root sub-system.
On the polynomial level, this property is reflected by van Diejen's
\cite{vD1} observation that the families of Macdonald polynomials
associated with classical root systems are special cases
or limit cases of the Koornwinder polynomials. In fact,
there is a large class of interesting
families of multivariable orthogonal polynomials associated with
classical root systems which are
degenerate cases of the Koornwinder polynomials, see e.g. \cite{KooS} and
\cite{vD0}. 
This supports the idea that the role of the Koornwinder polynomials
in the theory of multivariable orthogonal polynomials associated with
classical root systems is similar to the important and dominant role of the  
Askey-Wilson \cite{AW} polynomials 
in the theory of one variable
(basic) hypergeometric orthogonal polynomials.

In this paper we continue the affine Hecke algebra approach to the 
theory of Koornwinder polynomials. In particular,  Sahi's \cite{Sa2} 
bi-orthogonality relations for the non-symmetric Koornwinder polynomials
are extended to the case of continuous parameter values, and the corresponding
diagonal terms are evaluated explicitly.
We furthermore derive the evaluation formulas for the non-symmetric
Koornwinder polynomials. Anti-symmetric
Koornwinder polynomials are defined, 
and the explicit connection between
the non-symmetric and the (anti-)symmetric theory is established.
This leads to new derivations of Koornwinder's \cite{Ko},  
van Diejen's \cite{vD2} and Sahi's \cite{Sa}
results on the orthogonality relations, quadratic norm
evaluations and evaluation formulas for the symmetric Koornwinder
polynomials. We also shortly discuss 
the analogue of Weyl's character formula for Koornwinder polynomials,
and the (closely related) shift operators.
 
Instead of using shift operators to evaluate the diagonal terms for
the non-symmetric Koornwinder polynomials,
we use a method which is motivated by Cherednik's \cite{C4}
beautiful approach for proving Opdam's \cite{O} inversion formula of the
non-symmetric Harish-Chandra transform (the so-called Cherednik-Opdam
transform). This method amounts to evaluating the diagonal terms using
an explicit description of the action of the double affine Hecke
algebra on non-symmetric Koornwinder polynomials in terms of
operators acting on the spectral parameter.
This approach reveals interesting new structures, such as
an expression of the diagonal terms as residues of
the bi-orthogonality weight function. 
This method  
is also expected to be an important tool for obtaining
a better understanding
of $q$-analogues of the Cherednik-Opdam transform, see e.g. 
\cite{Cone} and \cite{KS} for some preliminary considerations 
in the rank one setting.

The results of this paper extend the results of Macdonald \cite{M3} and 
Cherednik \cite{C3} on the 
bi-orthogonality relations, diagonal term evaluations
and evaluation formulas for non-symmetric Macdonald polynomials
associated with classical reduced root systems, as well as the results of
Noumi \& Stokman \cite{NS}, in which the rank one setting was treated in
detail.

{\it Acknowledgements:} The author is supported by a fellowship from the Royal
Netherlands Academy of Arts and Sciences (KNAW). The research was done
during the author's stay at Universit{\'e} Pierre et Marie Curie (Paris VI) 
and  Institut de Recherche Math{\'e}matique Avanc{\'e}e (Strasbourg)
in France, supported by a NWO-TALENT
stipendium of the Netherlands Organization for Scientific Research
(NWO) and by
the EC TMR network 
``Algebraic Lie Representations'', grant no. ERB FMRX-CT97-0100.
The author thanks Masatoshi Noumi for fruitful discussions, and
Siddhartha Sahi and Akinori Nishino for drawing his attention to 
the papers \cite{Sa2} and \cite{NUKW}, respectively.


\section{The affine root system of type $C^\vee C_n$}

In this section we discuss the affine
root system of type $C^\vee C_n$, which was introduced by Macdonald
in \cite{M1}.
Let $n$ be a positive integer $\geq 2$.
Let $V=({\mathbb{R}}^n,\langle .,. \rangle)$ be
Euclidean $n$-space with orthonormal basis $\{\epsilon_i\}_{i=1}^n$.
We write $\widehat{V}$ for the affine linear transformations from $V$
to ${\mathbb{R}}$.
As a vector space, $\widehat{V}$ can be identified with $V\oplus
{\mathbb{R}}\delta$, where vectors in $V$ are considered as linear functionals
on $V$ via the scalar product $\langle .,. \rangle$, and where $\delta$
is the function identically equal to one on $V$. We extend the
scalar product $\langle .,. \rangle$ to a positive semi-definite
form on $\widehat{V}$  by requiring that the constant function
$\delta$ is in the radical of $\langle .,. \rangle$.

Let $S\subset \widehat{V}$ be the subset
\begin{equation}
\begin{split}
S=&\{\pm \epsilon_i+\frac{m}{2}\delta,\, \pm 2\epsilon_i+m\delta \,\, |
\,\,
m\in {\mathbb{Z}},\, i=1,\ldots,n\}\\
&\cup \{\pm \epsilon_i\pm \epsilon_j+m\delta \,\, | \,\, m\in
{\mathbb{Z}},\, 1\leq i<j\leq n\},
\end{split}
\end{equation}
where all the sign combinations occur. Let
${\mathcal{W}}={\mathcal{W}}(S)$
be the sub-group of $\hbox{GL}_{\mathbb{R}}(\widehat{V})$ generated by the
reflections $s_\beta$ ($\beta\in S$), where
\[s_f(g)= g-\langle g,f^\vee\rangle f,\qquad f\in
\widehat{V}\setminus {\mathbb{R}}\delta,\,\, g\in \widehat{V},
\]
and where $f^\vee=2f/\langle f,f\rangle$ is the co-root of $f$.
Observe that $s_f(g)=g\circ \widetilde{s}_f^{-1}$, with
$\widetilde{s}_f: V\rightarrow V$ the orthogonal reflection in
the affine hyperplane $f^{-1}(0)$.

By \cite{M1}, $S\subset \widehat{V}$ is an irreducible, affine root system.
In particular, $\langle \alpha,\beta^\vee\rangle\in {\mathbb{Z}}$
for all $\alpha,\beta\in S$, and 
$S$ is stable under the action of ${\mathcal{W}}$.
The sub-group ${\mathcal{W}}\subset \hbox{GL}(\widehat{V})$ is
called the affine Weyl group of $S$.

Let $R$ be the inmultiplyable
roots in $S$ and $R^\vee\subset S$ the corresponding co-root system.
Then $R$ and $R^\vee$ are irreducible, reduced
affine root systems in $\widehat{V}$, with affine Weyl group
${\mathcal{W}}$.
The projection $\Sigma\subset R$ of $R$ on $V$ along the direct
sum decomposition $\widehat{V}=V\oplus {\mathbb{R}}\delta$ is
an irreducible root system of type $C_n$ with Weyl group
$W=S_n\ltimes (\pm 1)^n\subset {\mathcal{W}}$ given by permutations and sign
changes of the fixed basis $\{\epsilon_i\}_{i=1}^n$ of $V$ (here $S_n$ 
denotes the symmetric group in $n$ letters). 
Due to the (non-disjoint) union $S=R\cup
R^\vee$ of $S$ into the 
reduced affine root sub-system $R$ of type $\widetilde{C}_n$
and its co-root system, we call $S$ of type $C^\vee C_n$, cf. \cite{M1}.

Let $Q^\vee$ be the co-root lattice of $\Sigma$, which coincides
with the weight lattice $\Lambda$ of $\Sigma$. In fact, $Q^\vee=\Lambda$ is
the full ${\mathbb{Z}}$-lattice in $V$ with basis $\{\epsilon_i\}_{i=1}^n$.
Then
\[{\mathcal{W}}=W\ltimes \tau(Q^\vee),
\]
where $\tau(v)\in \hbox{GL}(\widehat{V})$
($v\in V$) is the translation operator
defined by $\tau(v)f=f+\langle v,f\rangle\delta$ for $f\in \widehat{V}$.
Observe that $\tau(v)f=f\circ \widetilde{\tau}_v^{-1}$ with
$\widetilde{\tau}_v: V\rightarrow V$ 
given by $\widetilde{\tau}_v(\lambda)= \lambda-v$.

\begin{rem}\label{tau}
In Macdonald's \cite{M3} and Cherednik's \cite{C1}--\cite{C5}
work the translation operator
$\tau(v)$ in fact corresponds to $\tau(-v)=\tau(v)^{-1}$.
Later on, this change of convention (which is related to conjugation with
the largest Weyl group element $\sigma\in W$) causes
certain changes of signs compared with Cherednik's and
Macdonald's theory, see e.g. remark \ref{changestau}{\bf (ii)}.
\end{rem}

We fix a basis $\{a_i\}_{i=0}^n$ of $R$ by
\[a_0=\delta-2\epsilon_1,\,\, a_i=\epsilon_i-\epsilon_{i+1}\,\,
(i=1,\ldots,n-1),\,\,\, a_n=2\epsilon_n.
\]
Observe that
$\{a_0^\vee=a_0/2, a_1,\ldots,a_{n-1}, a_n^\vee=a_n/2\}$
is a basis of $R^\vee$, as well as of $S$. Furthermore,
$\{a_i\}_{i=1}^n$ is a basis of the gradient
root system $\Sigma$. We write $\Sigma^+$ (respectively
$\Sigma^-$) for the corresponding positive (respectively negative) roots in
$\Sigma$, and $\Lambda^+=\bigoplus_{i=1}^n{\mathbb{Z}}_+\omega_i$ for the
corresponding cone of dominant weights of $\Sigma$.
Here $\omega_i=\epsilon_1+\cdots
+\epsilon_i$ $(i=1,\ldots,n)$ are the fundamental weights of $\Lambda$,
i.e. $\langle \omega_i, a_j^\vee\rangle=\delta_{i,j}$ for all
$i,j=1,\ldots,n$, where $\delta_{i,j}$ is the Kronecker delta.
We furthermore write $Q^{\vee,+}$ for the positive span of the
simple co-roots $a_i^\vee$ ($i=1,\ldots,n$).

Let $R^+$ (respectively $R^-$) be the positive (respectively negative)
roots of $R$ with respect to the basis of the previous paragraph.
In particular, $R^+=\Sigma^+\cup
\{\beta\in R \, | \, \beta(0)>0\}$.

The affine Weyl group ${\mathcal{W}}$ is generated by the simple
reflections $s_i=s_{a_i}$ ($i=0,\ldots,n$). In fact, ${\mathcal{W}}$
is isomorphic to the Coxeter group with generators $s_i$
($i=0,\ldots,n$) satisfying $s_i^2=1$ and the braid relations
$s_is_{i+1}s_is_{i+1}=s_{i+1}s_is_{i+1}s_i$ ($i=0$, $i=n-1$),
$s_is_{i+1}s_i=s_{i+1}s_is_{i+1}$ ($i=1,\ldots, n-2$) and
$s_is_j=s_js_i$ for $|i-j|\geq 2$.

With our present conventions,
Lusztig's formula \cite[1.4(a)]{L}
for the length of an element in ${\mathcal{W}}$
is given by
\begin{equation}\label{length}
l(\tau(\lambda)w)=\sum_{\alpha\in\Sigma^+}|-\langle
\lambda,w\alpha\rangle+\chi(w\alpha)|,\qquad \lambda\in \Lambda, \,\,
w\in W,
\end{equation}
where $\chi(\alpha)=1$ if $\alpha\in \Sigma^-$ and $=0$
otherwise.

We write $\Sigma=\Sigma_m\cup \Sigma_l$ for the decomposition of $\Sigma$
into $W$-orbits, where $\Sigma_m$ (respectively $\Sigma_l$) is the
set of roots of length two (respectively four). We furthermore set
$\Sigma_s=\frac{1}{2}\Sigma_l$. There are five
${\mathcal{W}}$-orbits in $S$, namely
\begin{equation}\label{orbit}
\begin{split}
{\mathcal{W}}a_0^\vee&=\bigl(\frac{1}{2}+{\mathbb{Z}}\bigr)\delta+\Sigma_s,
\qquad
{\mathcal{W}}a_0=\bigl(1+2{\mathbb{Z}}\bigr)\delta+\Sigma_l,\\
{\mathcal{W}}a_i&={\mathbb{Z}}\delta+\Sigma_m\qquad (i\in \{1,\ldots,n-1\}
\,\, \hbox{ arbitrary}),\\
{\mathcal{W}}a_n^\vee&={\mathbb{Z}}\delta+\Sigma_s,\qquad
{\mathcal{W}}a_n=2{\mathbb{Z}}\delta+\Sigma_l.
\end{split}
\end{equation}
Observe that $R$ (respectively $R^\vee$) has three
${\mathcal{W}}$-orbits, namely ${\mathcal{W}}a_0$, ${\mathcal{W}}a_i$
and ${\mathcal{W}}a_n$ (respectively ${\mathcal{W}}a_0^\vee$, 
${\mathcal{W}}a_i$
and ${\mathcal{W}}a_n^\vee$), where $i\in\{1,\ldots,n-1\}$ is
arbitrary.

For later purposes, we define an action of ${\mathcal{W}}$
on $V$ which extends the canonical $W$-action on $V$. It suffices to
specify the action of the simple reflection $s_0$ on $V$, which we
take to be
\[s_0.x=(-1-x_1,x_2,\ldots,x_n),
\]
where $x_i=\langle x,\epsilon_i\rangle$. Observe that
$\Lambda\subset V$ is ${\mathcal{W}}$-stable,
and that $\tau(\lambda).x=x+\lambda$ for $\lambda\in \Lambda$.
We denote this action of ${\mathcal{W}}$ on $V$ with a dot and we
call it the {\it dot-action}, in order to avoid confusion with the
canonical action of ${\mathcal{W}}$ on $\widehat{V}$ and its induced
dual action on $V$.


\section{The (double) affine Hecke algebra}

Let ${\mathcal{A}}$ be the group algebra of the weight lattice $\Lambda$.
We write $x^{\lambda}$ ($\lambda\in \Lambda$) for the canonical
basis of ${\mathcal{A}}$, so that $x^0=1$ is the unit element in 
${\mathcal{A}}$
and $x^{\lambda}x^{\mu}=x^{\lambda+\mu}$
for all $\lambda,\mu\in \Lambda$.
The group algebra ${\mathcal{A}}$ is isomorphic to the Laurent
polynomials in the $n$
independent indeterminates $x_i=x^{\epsilon_i}$ ($i=1,\ldots,n$).
Let $q\in {\mathbb{C}}\setminus \{0\}$ be generic complex (in
particular, not a root of unity) and let $q^{1/2}$ be a fixed square root 
of $q$.
We write $x^{\mu+c\delta}=q^cx^{\mu}$
for $\mu\in \Lambda$ and $c\in \frac{1}{2}{\mathbb{Z}}$. Then the assignment
$w(x^{\mu})=x^{w\mu}$ for $w\in {\mathcal{W}}$ and
$\mu\in \Lambda\subset \widehat{V}$
extends by linearity to an action of ${\mathcal{W}}$ on
${\mathcal{A}}$. In particular,
the action of the simple reflections $s_i$ ($i=0,\ldots,n$) on ${\mathcal{A}}$
is given by
\begin{equation}\label{acA}
\begin{split}
(s_0f)(x)&=f(qx_1^{-1},x_2,\ldots,x_n),\\
(s_if)(x)&=f(x_1,\ldots,x_{i-1},x_{i+1},x_i,x_{i+2},\ldots,x_n)\qquad
(i=1,\ldots,n-1),\\
(s_nf)(x)&=f(x_1,\ldots,x_{n-1},x_n^{-1}),
\end{split}
\end{equation}
where $f\in {\mathcal{A}}$ and $x=(x_1,\ldots,x_n)$.
In particular, the translation operators $\tau(\mu)$
($\mu\in \Lambda$) act as $q$-difference operators:
$\tau(\mu)(x^{\lambda})=q^{\langle \mu, \lambda\rangle}x^{\lambda}$
for all $\lambda,\mu\in \Lambda$.

The Noumi \cite{N} representation is a five($=\#\{{\mathcal{W}}-\hbox{orbits
in } S\}$) parameter deformation of the above action of ${\mathcal{W}}$ on
the group algebra ${\mathcal{A}}$.
We incorporate the five extra degrees of freedom in a
so-called multiplicity function
$\mathbf{t}=(t_\beta)_{\beta\in S}$ of $S$, which is a
${\mathcal{W}}$-invariant map from
$S$ to ${\mathbb{C}}\setminus \{0\}$ (so $t_{w\beta}=t_\beta$
for all $w\in {\mathcal{W}}$ and all $\beta\in
S$). We furthermore set $t_f=1$ if $f\in \widehat{V}\setminus S$.
A multiplicity function is thus uniquely determined by the
five values $t_{a_0^\vee}$, $t_{a_0}$, $t_{a_i}$ ($i\in \{1,\ldots,n-1\}$
arbitrary), $t_{a_n}$ and $t_{a_n^\vee}$. In order to avoid
cumbersome notations, we sometimes write $t_{i}$ (respectively $t_i^\vee$)
for $t_{a_i}$ (respectively $t_{a_i^\vee}$) and we write $t$
for the value of $t_j=t_j^\vee$ with $j\in \{1,\ldots,n-1\}$.
Furthermore, we use the short-hand 
notation $\mathbf{k}=(t_\beta)_{\beta\in R}\simeq
(t_{0},t,t_{n})$ and $\mathbf{k}^\vee=
(t_\beta)_{\beta\in R^\vee}\simeq (t_{0}^\vee,
t,t_{n}^\vee)$
for the corresponding multiplicity functions of $R$ and $R^\vee$,
respectively. We assume throughout the paper that the values of the
multiplicity function $\mathbf{t}$ are generically complex.

In the Noumi representation, the role of the affine Weyl group ${\mathcal{W}}$
is replaced by the affine Hecke algebra of type $\widetilde{C}_n$,
which is defined as follows.

\begin{Def} The affine Hecke algebra 
$H=H(R;\mathbf{k})$
of type $\widetilde{C}_n$ is the unital, associative
algebra with generators $T_0,\ldots,T_n$ and relations
\begin{equation}
\label{relfund1}
(T_i-t_{i})(T_i+t_{i}^{-1})=0,\qquad (i=0,\ldots,n),
\end{equation}
and the braid relations
\begin{equation}\label{relfund2}
\begin{split}
T_iT_{i+1}T_iT_{i+1}&=T_{i+1}T_iT_{i+1}T_i,\qquad (i=0,i=n-1),\\
T_iT_{i+1}T_i&=T_{i+1}T_iT_{i+1},\qquad (i=1,\ldots,n-2),\\
T_iT_j&=T_jT_i,\qquad |i-j|\geq 2.
\end{split}
\end{equation}
\end{Def}
We call \eqref{relfund1} and \eqref{relfund2} the $H(R;\mathbf{k})$-relations
for the $(n+1)$-tuple $(T_0,\ldots,T_n)$. Furthermore, 
we write $H(R^\vee;\mathbf{k}^\vee)$ for the affine Hecke algebra $H$
in which the parameter $t_i$ is replaced by $t_i^{\vee}$ for
$i=0,\ldots,n$.

We recall here some of the basic properties of the affine Hecke algebra $H$,
see Lusztig \cite{L}
for details and for a general discussion on affine Hecke algebras.

For a reduced expression $w=s_{i_1}\cdots s_{i_r}$ of $w\in {\mathcal{W}}$
we set $T_w=T_{i_1}\cdots T_{i_r}$.
This is independent of the choice of reduced
expression by the braid relations \eqref{relfund2}
for the $T_i$, and $\{T_w\}_{w\in {\mathcal{W}}}$
is a linear basis of $H$.

For $\lambda\in \Lambda^+$ we set
$Y^{\lambda}=T_{\tau(\lambda)}$, and for $\lambda=\mu-\nu\in \Lambda$
with $\mu,\nu\in \Lambda^+$ we set $Y^{\lambda}=Y^{\mu}(Y^{\mu})^{-1}$.
The length identity \eqref{length} implies that the
$Y^{\lambda}$ ($\lambda\in \Lambda$)
are well-defined (i.e. independent of the choice of
decomposition $\lambda=\mu-\nu$). Furthermore, the
sub-space ${\mathcal{A}}_Y=
\hbox{span}\{Y^{\lambda} \, | \, \lambda\in \Lambda\}$
is a commutative subalgebra of $H$ isomorphic to ${\mathcal{A}}$
(in particular, $Y^0=1$ and $Y^{\lambda}Y^{\mu}=Y^{\lambda+\mu}$
for all $\lambda,\mu\in \Lambda$). We identify
$f(x)=\sum_{\lambda}c_{\lambda}x^{\lambda}\in {\mathcal{A}}$
with $f(Y)=\sum_{\lambda}c_{\lambda}Y^{\lambda}\in {\mathcal{A}}_Y$
in the remainder of the paper. 
We write $Y_i=Y^{\epsilon_i}$, which
corresponds with $x_i$ under the identification of ${\mathcal{A}}_Y$
with $\mathcal{A}$.

The Hecke algebra $H_0=H_0(\Sigma;t,t_n)$ of the finite Weyl
group $W$ can be identified
with the subalgebra of $H$ generated by $T_i$ ($i=1,\ldots,n$).
Then $\{T_w\}_{w\in W}$ is a linear basis of $H_0$ and
\[H\simeq H_0\otimes {\mathcal{A}}_Y\simeq
{\mathcal{A}}_Y\otimes H_0
\]
as vector spaces by multiplication. The commutation relations
between $T_i\in H_0$ ($i=1,\ldots,n$) and $f(Y)\in {\mathcal{A}}_Y$
are given by the formulas
\begin{equation}\label{lusztig}
\begin{split}
T_if(Y)-(s_if)(Y)T_i&=(t-t^{-1})
\left(\frac{f(Y)-(s_if)(Y)}{1-Y^{-a_i}}\right),\\
T_nf(Y)-(s_nf)(Y)T_n&=\bigl((t_{n}-t_{n}^{-1})+
(t_{0}-t_{0}^{-1})Y_n^{-1}\bigr)
\left(\frac{f(Y)-(s_nf)(Y)}
{1-Y_n^{-2}}\right)
\end{split}
\end{equation}
for $i=1,\ldots,n-1$, see \cite[prop. 3.6]{L}.
These commutation relations can be used to prove inductively that
\begin{equation}\label{YT}
Y_i=T_i\cdots T_{n-1}T_nT_{n-1}\cdots T_1T_0T_1^{-1}T_2^{-1}\cdots
T_{i-1}^{-1},\qquad i\in\{1,\ldots,n\}
\end{equation}
in $H$, see Noumi \cite{N} or Sahi \cite[(11)]{Sa}.
The expression \eqref{YT} is the analogue in $H$ of the reduced
expression
\begin{equation}
\label{taui}
\tau(\epsilon_i)=s_i\cdots s_{n-1}s_ns_{n-1}\cdots
s_1s_0s_1\cdots s_{i-1},\qquad i\in\{1,\ldots,n\}
\end{equation}
in ${\mathcal{W}}$. We define a rational function
$v_\beta(x)=v_\beta(x;\mathbf{t};q)\in
{\mathbb{C}}(x)=\hbox{Quot}({\mathcal{A}})$ by
\begin{equation}\label{v}
v_\beta(x;\mathbf{t};q)=\frac{\bigl(1-t_{\beta}t_{\beta/2}\,x^{\beta/2}\bigr)
\bigl(1+t_{\beta}t_{\beta/2}^{-1}x^{\beta/2}\bigr)}{\bigl(1-x^\beta\bigr)},
\qquad
\beta\in R.
\end{equation}
Observe that for $\beta\in R$ with $\beta/2\not\in S$, the expression
\eqref{v} reduces to $v_\beta(x)=
(1-t_\beta^2x^\beta)/(1-x^\beta)$ since $t_{\beta/2}=1$.
The following crucial theorem was proved by Noumi \cite{N}.
\begin{thm}[The Noumi representation]\label{noumi}
The assignment
\begin{equation*}
T_i\mapsto
t_{i}+t_{i}^{-1}v_{a_i}(x;\mathbf{t};q)\bigl(s_i-\hbox{id}\bigr)
\in \hbox{End}_{\mathbb{C}}({\mathcal{A}})
\end{equation*}
for $i=0,\ldots,n$ uniquely extends to a representation
$\pi_{\mathbf{t},q}: H(R;\mathbf{k})\rightarrow
\hbox{End}_{\mathbb{C}}({\mathcal{A}})$.
\end{thm}
The commutation relations \eqref{lusztig} play a crucial role in
the proof of theorem \ref{noumi}, compare with the argument in \cite{C1}
and \cite[(4.6)]{M3} in case of reduced root systems.

We write $T_i$ for the image of $T_i\in H(R;\mathbf{k})$
under the Noumi representation $\pi_{\mathbf{t},q}$ for $i=0,\ldots,n$
if no confusion is possible, and we 
call them the difference-reflection operators
associated with $S$.

We are now in a position to recall Sahi's \cite{Sa} definition of the double 
affine Hecke algebra. 

\begin{Def}
The {\it double affine Hecke algebra}
${\mathcal{H}}={\mathcal{H}}(S;\mathbf{t};q)$
is the sub-algebra of $\hbox{End}_{\mathbb{C}}({\mathcal{A}})$
generated by $\pi_{\mathbf{t},q}(H(R;\mathbf{k}))$ and
${\mathcal{A}}$, where the elements in ${\mathcal{A}}$ are considered as 
multiplication operators in 
$\hbox{End}_{\mathbb{C}}({\mathcal{A}})$.
\end{Def}

We end this section by giving an alternative presentation of
${\mathcal{H}}$ (also different from Sahi's \cite[sect. 3]{Sa} presentation), 
which emphasizes its close connection with the affine root system $S$.

We write $f(z)=\sum_{\lambda}c_{\lambda}z^{\lambda}\in 
\hbox{End}_{\mathbb{C}}({\mathcal{A}})$ for
the multiplication operator associated with
the Laurent polynomial 
$f(x)=\sum_{\lambda}c_{\lambda}x^{\lambda}\in {\mathcal{A}}$.
In particular, $z^{\lambda+m\delta}=q^mz^{\lambda}$ ($\lambda\in
\Lambda$, $m\in \frac{1}{2}{\mathbb{Z}}$) is the multiplication operator
associated with $x^{\lambda+m\delta}=q^mx^{\lambda}\in {\mathcal{A}}$.
Then in ${\mathcal{H}}$ we have the
commutation relations
\begin{equation}\label{extra}
f(z)T_i-T_i(s_if)(z)=
\frac{(t_{a_i}-t_{a_i}^{-1})+(t_{a_i/2}-t_{a_i/2}^{-1})z^{a_i/2}}
{1-z^{a_i}}\bigl(f(z)-(s_if)(z)\bigr)
\end{equation}
for $i=0,\ldots,n$ and $f\in {\mathcal{A}}$.
This follows from the fact that
the difference-reflection operator $T_i$ can be rewritten as
\[
T_i=
t_{a_i}s_i+\frac{(t_{a_i}-t_{a_i}^{-1})+
(t_{a_i/2}-t_{a_i/2}^{-1})x^{a_i/2}}
{1-x^{a_i}}(\hbox{id}-s_i)
\]
for $i=0,\ldots,n$.
\begin{thm}\label{double}
The double affine Hecke algebra ${\mathcal{H}}(S;\mathbf{t};q)$
is isomorphic to the unital, associative algebra ${\mathcal{F}}(\mathbf{t};q)$
with generators $V_0^\vee, V_0$, $V_i$
\textup{(}$i=1,\ldots,n$\textup{)} and $V_n^\vee$, satisfying
\begin{enumerate}
\item[{1.}] The $H(R;\mathbf{k})$-relations for
$(V_0,V_1,\ldots,V_{n-1},V_n)$.
\item[{2.}] The $H(R^\vee;\mathbf{k}^\vee)$-relations for
$(V_0^\vee,V_1,\ldots,V_{n-1},V_n^\vee)$.
\item[{3.}] \textup{(}Compatibility conditions\textup{)}.
$V_n^\vee V_nV_{n-1}\cdots V_1
V_0V_0^\vee V_1V_2\cdots V_{n-1}=q^{-1/2}$ and
$\lbrack V_0,V_n^\vee\rbrack=0=\lbrack V_0^\vee,V_n\rbrack$. 
\end{enumerate}
The algebra isomorphism
$\phi: {\mathcal{F}}(\mathbf{t};q)\rightarrow {\mathcal{H}}(S;\mathbf{t};q)$
is explicitly given by
$\phi(V_i)=T_i$ \textup{(}$i=0,\ldots,n$\textup{)},
$\phi(V_0^\vee)=T_0^{-1}z^{-a_0^\vee}=
q^{-1/2}T_0^{-1}z_1$ and
$\phi(V_n^\vee)=z^{-a_n^\vee}T_n^{-1}=z_n^{-1}T_n^{-1}$.
\end{thm}
\begin{proof}
For the existence of $\phi$, we need to check that
the elements $T_0^\vee=T_0^{-1}z^{-a_0^\vee}$, $T_j$
($j=0,\ldots,n$) and $T_n^\vee=z^{-a_n^\vee}T_n^{-1}$ in ${\mathcal{H}}$
respect the defining relations of the generators $V_0^\vee$, $V_j$
($j=0,\ldots,n$), $V_n^\vee$ in ${\mathcal{F}}$.
The $H(R;\mathbf{k})$-relations for $(T_0,\ldots,T_n)$ is precisely the
content of theorem \ref{noumi}.
The $H(R^\vee; \mathbf{k}^\vee)$-relations for the $(n+1)$-tuple
$(T_0^\vee, T_1,\ldots,T_{n-1}, T_n^\vee)$ follows easily from
\eqref{extra} (see also \cite[sect. 3]{Sa}).
By \eqref{extra} it follows inductively that
\begin{equation*}
\begin{split}
z_i&=T_i^{-1}T_{i+1}^{-1}\cdots
T_{n}^{-1}(T_n^\vee)^{-1}T_{n-1}^{-1}\cdots T_i^{-1}\\
&=q^{1/2}T_{i-1}\cdots T_1T_0T_0^\vee T_1\cdots T_{i-1}
\end{split}
\end{equation*}
for $i=1,\ldots,n$ in ${\mathcal{H}}$. In particular,
the identity $z_n^{-1}z_n=1$
in ${\mathcal{H}}$ shows that 
the generators $T_0^\vee$, $T_j$ ($j=0,\ldots,n$) and $T_n^\vee$
satisfy the compatibility condition in ${\mathcal{H}}$. Hence the algebra
homomorphism $\phi$ exists.

We write $w_i=q^{1/2}V_{i-1}\cdots 
V_1V_0V_0^\vee V_1\cdots V_{i-1}\in {\mathcal{F}}$
for $i\in\{1,\ldots,n\}$, so that $\phi(w_i)=z_i$ for $i=1,\ldots,n$.
Since the $z_i$ ($i=1,\ldots,n$) and the
$T_j=\phi(V_j)$ ($j=0,\ldots,n$) generate ${\mathcal{H}}$
as an algebra, we see that $\phi$ is surjective.
On the other hand, all fundamental relations of ${\mathcal{H}}(S;\mathbf{t};q)$
as given by Sahi \cite[sect. 3]{S} can be easily checked for
the generators $w_i$ ($i=1,\ldots,n$) and $V_j$ ($j=0,\ldots,n$)
of ${\mathcal{F}}$. This implies the injectivity of
$\phi$.
\end{proof}

\begin{rem}\label{Spresent}
The presentation of ${\mathcal{H}}(S;\mathbf{t};q)$
as given in theorem \ref{double} clearly reflects
the structure of the
underlying non-reduced affine root system $S$.  In particular, the
first part of the compatibility condition can be recovered 
from the root data as follows.
We put the simple roots for the indivisible roots $R^\vee$
(respectively for the inmultiplyable roots $R$)
above (respectively below) the corresponding vertices of the extended
Dynkin diagram:
\begin{equation*}
\begin{CD}
\underset{\displaystyle{a_0}}{\stackrel{\displaystyle{a_0^\vee}}
{\displaystyle{\circ}}} @=
\underset{\displaystyle{a_1}}{\stackrel{\displaystyle{a_1}}{\displaystyle{\circ}}}
\frac{\qquad\quad}{\qquad\quad}
\underset{\displaystyle{a_2}}{\stackrel{\displaystyle{a_2}}{\displaystyle{\circ}}}
\frac{\qquad\quad}{\qquad\quad}\,\cdots\cdots\,\frac{\qquad\quad}{\qquad\quad}
\underset{\displaystyle{a_{n-1}}}{\stackrel{\displaystyle{a_{n-1}}}
{\displaystyle{\circ}}}
@= \underset{\displaystyle{a_n}}{\stackrel{\displaystyle{a_n^\vee}}
{\displaystyle{\circ}}}
\end{CD}
\end{equation*}
Now we attach $V_i(=T_i)$ to the simple roots $a_i$,
$V_0^\vee(=T_0^\vee=T_0^{-1}z^{-a_0^\vee})$ to the co-root $a_0^\vee$
and $V_n^\vee(=T_n^\vee=z^{-a_n^\vee}T_n^{-1})$ to the co-root
$a_n^\vee$ in the above diagram. We call them the {\it simple
generators} of ${\mathcal{F}}\simeq {\mathcal{H}}$.
Then
the compatibility condition amounts to the following
rule: {\it multiplying simple generators in the order of appearance
of a single walk around the diagram in clockwise direction, gives} $q^{-1/2}$.
The point of departure for the walk is
irrelevant, since the compatibility condition is equivalent to the
compatibility condition in which the factors in its left-hand side
are permuted cyclically.
\end{rem}


\section{Non-symmetric Koornwinder polynomials and triangularity}

In this section we show that
the $Y$-operators $Y^{\lambda}$ ($\lambda\in \Lambda$)
act as triangular operators
under the Noumi representation $\pi_{\mathbf{t},q}$. We use this
triangularity property to redefine Sahi's \cite{Sa} non-symmetric Koornwinder
polynomials. The advantage of this method is that triangularity
properties of the non-symmetric Koornwinder polynomials are
automatically incorporated in their definition, in contrast with Sahi's \cite{Sa},
\cite{Sa2} approach.

Let $\lambda^+\in \Lambda^+$ for $\lambda\in\Lambda$
be the unique dominant weight in the
orbit $W\lambda$. We will be needing the following two partial orders
on the weight lattice $\Lambda$.
\begin{Def} Let $\lambda,\mu\in \Lambda$.\\
{\bf (i)} We write $\lambda\leq\mu$ if
$\mu-\lambda\in Q^{\vee,+}$
\textup{(}and $\lambda<\mu$
if $\lambda\leq\mu$ and $\lambda\not=\mu$\textup{)}.\\
{\bf (ii)} We write $\lambda\preceq\mu$
if $\lambda^+<\mu^+$, or if $\lambda^+=\mu^+$ and $\lambda\leq\mu$
\textup{(}and $\lambda\prec\mu$ if $\lambda\preceq\mu$ and
$\lambda\not=\mu$\textup{)}.
\end{Def}

\begin{lem}\label{technical}
Let $\mu\in \Lambda$ and $\alpha\in\Sigma^+$.

If $\langle \mu,\alpha\rangle \geq 2$, then
$\mu-r\alpha^\vee\prec \mu$ for
$r=1,\ldots, \langle \mu,\alpha\rangle-1$.

If $\langle \mu,\alpha\rangle\leq -2$, then $\mu+r\alpha^\vee\prec \mu$
for $r=1,\ldots,-\langle \mu,\alpha\rangle-1$.
\end{lem}
\begin{proof}
We write $m_{\alpha}=\langle \mu,\alpha\rangle\in {\mathbb{Z}}$.
Suppose that $m_\alpha\geq 2$ and write $\mu_r=\mu-r\alpha^\vee$
with $r\in \{1,\ldots,m_\alpha-1\}$. We show that $\mu_r^+<\mu^+$.

Let $w\in W$
such that $\mu_r^+=w\mu_r$.
If $w\alpha^\vee\in Q^{\vee,+}$,
then $\mu_r^+=w\mu-rw\alpha^\vee<w\mu\leq\mu^+$. On the other
hand, if $w\alpha^\vee\in -Q^{\vee,+}$,
then $\mu_r^+=w\mu-rw\alpha^\vee<w\mu-m_\alpha w\alpha^\vee
=(ws_\alpha)\mu\leq \mu^+$.
This proves the assertion for $m_\alpha\geq 2$.
The case $m_\alpha\leq -2$ can be obtained by applying
the previous case to $s_\alpha\mu$.
\end{proof}

For $\beta\in R$ we define
\begin{equation}
\label{R}
{\mathcal{R}}(\beta)=t_{\beta}s_\beta+
t_\beta^{-1}v_\beta(x)\bigl(1-s_\beta\bigr)\in
\hbox{End}_{\mathbb{C}}({\mathcal{A}}),
\end{equation}
where $v_\beta(\cdot)$ is given by \eqref{v}.
Let $\epsilon: {\mathbb{Z}}\rightarrow \{\pm 1\}$ be the function
which maps a positive integer to $1$ and a strictly negative
integer to $-1$.
\begin{lem}\label{triangularR}
Let $\lambda\in \Lambda$. For 
$\beta=\alpha+m\delta\in R^+$ with $\alpha\in \Sigma^+$
we have
\[
{\mathcal{R}}(\beta)(x^{\lambda})=t_\beta^{\epsilon(\langle
\lambda,\beta\rangle)}x^{\lambda}+
\sum_{\mu\prec\lambda}c_{\lambda,\mu}x^{\mu}
\]
for certain constants $c_{\lambda,\mu}\in {\mathbb{C}}$.
\end{lem}
\begin{proof}
Let $D_\beta\in \hbox{End}_{\mathbb{C}}({\mathcal{A}})$
for $\beta\in R$ be the divided difference-reflection operator defined by
\[ D_{\beta}f=\frac{f-s_{\beta}f}{1-x^\beta},\qquad f\in {\mathcal{A}}.
\]
Then for all $\lambda\in \Lambda$ we have
\begin{equation*}
D_\beta(x^{\lambda})=
\begin{cases}
-x^{\lambda-\beta}-x^{\lambda-2\beta}-\cdots -x^{\lambda-\langle
\lambda,\beta^\vee\rangle \beta}\qquad &\hbox{if }\,\, \langle
\lambda,\beta^\vee\rangle>0,\\
0 &\hbox{if }\,\, \langle \lambda, \beta^\vee\rangle=0,\\
x^{\lambda}+x^{\lambda+\beta}+\cdots +
x^{\lambda-(1+\langle \lambda,\beta^\vee\rangle)\beta}
&\hbox{if }\,\, \langle \lambda,\beta^\vee\rangle<0.
\end{cases}
\end{equation*}
The proof follows now easily from lemma \ref{technical} and from the
definition \eqref{R} of ${\mathcal{R}}(\beta)$.
\end{proof}
Observe that ${\mathcal{R}}(a_i)=T_is_i$ for $i=0,\ldots,n$ and
that ${\mathcal{R}}(w(\beta))=w{\mathcal{R}}(\beta)w^{-1}$ for all
$w\in {\mathcal{W}}$
and all $\beta\in R$. Combined with \eqref{YT} and \eqref{taui}, we obtain
\begin{equation}\label{Yalternatief}
\begin{split}
Y_i=&{\mathcal{R}}(\epsilon_i-\epsilon_{i+1})
{\mathcal{R}}(\epsilon_i-\epsilon_{i+2})\cdots
{\mathcal{R}}(\epsilon_i-\epsilon_n){\mathcal{R}}(2\epsilon_i)\\
&\cdot{\mathcal{R}}(\epsilon_i+\epsilon_n)\cdots
{\mathcal{R}}(\epsilon_i+\epsilon_{i+1})
{\mathcal{R}}(\epsilon_i+\epsilon_{i-1})\cdots
{\mathcal{R}}(\epsilon_i+\epsilon_1)\\
&\cdot{\mathcal{R}}(\delta+2\epsilon_i)\tau(\epsilon_i)
{\mathcal{R}}(\epsilon_1-\epsilon_i)^{-1}\cdots
{\mathcal{R}}(\epsilon_{i-1}-\epsilon_i)^{-1}
\end{split}
\end{equation}
for $i=1,\ldots,n$, cf. \cite{N}.
Hence the triangularity of the factors ${\mathcal{R}}(\cdot)$ in
\eqref{Yalternatief} (see lemma \ref{triangularR})
implies the triangularity of $Y_i$
for $i=1,\ldots,n$, and hence of $Y^{\lambda}$ for all $\lambda\in \Lambda$.
For the explicit description of the diagonal
terms of the $Y$-operators, we need to introduce some
additional notations first.  

We write $f(y)$ for the value of $f\in {\mathcal{A}}$
at $y=(y_1,\ldots,y_n)\in ({\mathbb{C}}\setminus 
\{0\})^n$. In particular, $(y)^{\lambda+m\delta}=q^m(y)^{\lambda}$
($m\in \frac{1}{2}{\mathbb{Z}}$, $\lambda\in\Lambda$) is the value of
$x^{\lambda+m\delta}=q^mx^{\lambda}\in {\mathcal{A}}$ at $y$.
Conversely, we let $c^{\lambda}\in ({\mathbb{C}}\setminus \{0\})^n$
for $c\in {\mathbb{C}}\setminus \{0\}$ and $\lambda\in\Lambda$
be the vector $c^{\lambda}=(c^{\lambda_1},\ldots, c^{\lambda_n})$, where
$\lambda_i=\langle \lambda,\epsilon_i\rangle$. 

\begin{rem}
The brackets in the notation $(y)^{\lambda+m\delta}$ for the value of 
$x^{\lambda+m\delta}$ at $y\in ({\mathbb{C}}\setminus\{0\})^n$ will 
occasionally be omitted when $m=0$. For $m\not=0$ the brackets are needed
to distinguish the value of $x^{\lambda+m\delta}$
at $y^{-1}=(y_1^{-1},\ldots,y_n^{-1})$
from the value of $x^{-\lambda-m\delta}$ at $y$. 
\end{rem}

Let now
$\Sigma_m^+$ (respectively $\Sigma_l^+$) be the positive roots in
$\Sigma$ of squared length $2$ (respectively $4$). Let
\[\rho_m=2\sum_{i=1}^n(n-i)\epsilon_i,\qquad
\rho_l=\sum_{i=1}^n\epsilon_i
\]
be the sum of co-roots $\alpha^\vee$ with $\alpha\in\Sigma_m^+$
and $\alpha\in\Sigma_l^+$, respectively.
Let $\lambda\in \Lambda$ and set
\[\rho_m(\lambda)=\sum_{\alpha\in\Sigma_m^+}\epsilon(\langle
\lambda,\alpha\rangle)\alpha^\vee,\qquad
\rho_l(\lambda)=\sum_{\alpha\in\Sigma_l^+}\epsilon(\langle
\lambda,\alpha\rangle)\alpha^\vee.
\]
Then $\rho_m(\lambda)=\rho_m$ and $\rho_l(\lambda)=\rho_l$
for all dominant weights $\lambda\in \Lambda^+$.
We define now $\gamma_{\lambda}=\gamma_{\lambda}(\mathbf{k},q)\in
{\mathbb{C}}^n$ by
\begin{equation}\label{gl}
\gamma_{\lambda}=t_{0}^{\rho_l(\lambda)}.t_{n}^{\rho_l(\lambda)}
.t^{\rho_m(\lambda)}.q^{\lambda},\qquad
\lambda\in \Lambda,
\end{equation}
where the product is the usual dot-product in ${\mathbb{C}}^n$.
In particular,
\begin{equation}\label{positive}
\gamma_{\lambda}=(t_{0}t_{n}t^{2(n-1)}q^{\lambda_1},
t_{0}t_{n}t^{2(n-2)}q^{\lambda_2},\ldots,
t_{0}t_{n}q^{\lambda_n}),
\qquad \lambda\in \Lambda^+.
\end{equation}
Then lemma \ref{triangularR}, \eqref{Yalternatief} and
\eqref{orbit} lead to the following proposition.

\begin{prop}\label{triangularY}
For $\lambda\in \Lambda$ and $f(Y)\in {\mathcal{A}}_Y$, we have
\[
f(Y)(x^{\lambda})=f(\gamma_{\lambda})x^{\lambda}+\sum_{\mu\prec\lambda}
c_{\lambda,\mu}x^{\mu}
\]
for certain constants $c_{\lambda,\mu}$.
\end{prop}
We give now some
properties of the diagonal terms $\gamma_{\lambda}$ ($\lambda\in \Lambda$)
which will be used frequently in the remainder of the paper.
First of all, the diagonal terms of the $Y$-operators can be related to
the spectrum of the $Y$-operators as described in 
\cite[def. 2.5]{Sa} by observing that
\begin{equation}\label{linkdiag}
\rho_m(\lambda)=w_{\lambda}\rho_m,\qquad \rho_l(\lambda)=
w_{\lambda}\rho_l,\qquad \lambda\in \Lambda,
\end{equation}
where $w_{\lambda}\in W$ is the unique element of minimal length such
that $\lambda^+=w_{\lambda}^{-1}\lambda$, see \cite[prop. 2.10]{O}.

Secondly, the action of ${\mathcal{W}}$ on the diagonal terms
$\gamma_{\lambda}$ ($\lambda\in\Lambda$), induced from the dot-action of 
${\mathcal{W}}$ on $\Lambda$, is compatible
with the action of ${\mathcal{W}}$
on ${\mathcal{A}}$ in the following way.
We refer to \cite[thm. 5.3]{Sa} for the proof.
\begin{lem}\label{compatibleactions}
Let $f\in {\mathcal{A}}$. Let $\lambda\in \Lambda$ and
$i\in \{0,1,\ldots,n\}$ such that $s_i.\lambda\not=\lambda$.
Then $f(\gamma_{s_i.\lambda}^{-1})=(s_if)(\gamma_{\lambda}^{-1})$.
If furthermore $i\geq 1$, then also
$f(\gamma_{s_i.\lambda})=f(\gamma_{s_i\lambda})=(s_if)(\gamma_{\lambda})$.
\end{lem}

\begin{rem}\label{casestable}
Observe that the condition $s_i.\lambda\not=\lambda$ in 
lemma \ref{compatibleactions}
is always met for $i=0$. Let now $i\in \{1,\ldots,n\}$ and 
$\lambda\in \Lambda$
with $s_i.\lambda(=s_i\lambda)=\lambda$.
Then we have $\gamma_{\lambda}^{a_i}=\gamma_0^{a_i}$
($=t^2$ if $i<n$ and $=t_0^2t_n^2$ if $i=n$).
Hence $(s_if)(\gamma_{\lambda}^{\pm 1})=f({\gamma}_{\lambda,i}^{\pm 1})$
for $f\in {\mathcal{A}}$, with ${\gamma}_{\lambda,i}=
\gamma_{\lambda}.(\gamma_0^{-a_i})^{a_i^\vee}$. 
Observe that ${\gamma}_{\lambda,i}\not=
\gamma_{\mu}$ for all $\mu\in \Lambda$ by the generic conditions on the
parameters.
\end{rem}

By lemma \ref{compatibleactions} we have $f(\gamma_{\lambda})=
\bigl(w_{\lambda}^{-1}f\bigr)(\gamma_{\lambda^+})$ for
$f\in {\mathcal{A}}$ and $\lambda\in \Lambda$.
Combined with \eqref{positive} we conclude
that the diagonal terms $\gamma_{\lambda}$ ($\lambda\in \Lambda$) are
mutually different for generic values of $q$ and $\mathbf{k}$.
This leads to the following main result of this section.
\begin{thm}\label{nonsymmetric}
There exists a unique basis $\{P_{\lambda}\}_{\lambda\in \Lambda}$
of ${\mathcal{A}}$ such that
\begin{enumerate}
\item[{--}] $P_{\lambda}(x)=x^{\lambda}+\sum_{\mu\prec\lambda}
c_{\lambda,\mu}x^{\mu}$ for certain constants $c_{\lambda,\mu}$,
\item[{--}] $f(Y)P_{\lambda}=f(\gamma_{\lambda})P_{\lambda}$ for
all $f(Y)\in {\mathcal{A}}_Y$,
\end{enumerate}
for all $\lambda\in \Lambda$.
\end{thm}
\begin{Def}\label{defnonsym}
The Laurent polynomial $P_{\lambda}(\cdot)=P_{\lambda}(\cdot;\mathbf{t};q)$
\textup{(}$\lambda\in \Lambda$\textup{)}
is called the monic, non-symmetric Koornwinder polynomial
of degree $\lambda$.
\end{Def}
The terminology introduced in definition \ref{defnonsym}
stems from the close connection between the
Laurent polynomials $P_{\lambda}$ ($\lambda\in\Lambda$) 
and Koornwinder's \cite{Ko}
multivariable analogues of the Askey-Wilson polynomials, see
\cite{N} and \cite{Sa}, as well as section 5 and section 6.

\begin{rem}\label{changestau}
{\bf (i)} The second property of theorem \ref{nonsymmetric} already
characterizes the non-symmetric Koornwinder polynomial $P_{\lambda}$
up to a constant.
This characterizing property was used by Sahi \cite[def. 6.1]{Sa}
to introduce the non-symmetric Koornwinder polynomials. The
triangularity of the non-symmetric Koornwinder polynomials was
derived by Sahi \cite[sect. 6]{Sa2} using recursion formulas.

{\bf (ii)} If one uses Macdonald's \cite{M3} and
Cherednik's \cite{C1}--\cite{C5} convention for the translation 
operator $\tau$ (see
remark \ref{tau}), then the role of $Y^{\lambda}=T_{\tau(\lambda)}$
is
taken over by $T_{\tau(\sigma\lambda)}=T_{\sigma}^{-1}Y^{\lambda}T_{\sigma}$
for $\lambda\in \Lambda^+$, where
$\sigma\in W$ is the longest Weyl group element. The common
eigenfunctions then become $T_{\sigma}^{-1}P_{\lambda}$ ($\lambda\in
\Lambda$). {}From lemma \ref{triangularR} and theorem
\ref{nonsymmetric} it follows that the $T_{\sigma}^{-1}P_{\lambda}$
are triangular with respect to the partial
order on $\Lambda$ in which $\mu$ is less than $\nu$ iff 
$\sigma\mu\prec \sigma\nu$
(i.e. the anti-dominant weight is highest in each
$W$-orbit). This is in accordance with the triangular structure of
non-symmetric Macdonald polynomials, see \cite{M3} and \cite{C3}.

\end{rem}


\section{Symmetric Koornwinder polynomials}

In this section we recall Noumi's \cite{N} results on the affine
Hecke algebraic characterization of Koornwinder's \cite{Ko} multivariable
analogues of the Askey-Wilson polynomials. We present Noumi's results here in a
different order by making use of the triangularity
of the $Y$-operators, see proposition \ref{triangularY}.

We write ${\mathcal{A}}^W=
\{ f\in {\mathcal{A}} \, | \, wf=f\,\,\,\,\forall w\in
W\}$, and similarly ${\mathcal{A}}_Y^W$, where the action
is given by $w(Y^{\lambda})=Y^{w\lambda}$ for $w\in W$ and $\lambda\in
\Lambda$. A linear basis of ${\mathcal{A}}^W$
and ${\mathcal{A}}_Y^W$ is given by the
monomials $m_{\lambda}(x)=\sum_{\mu\in W\lambda}x^{\mu}$, respectively
$m_{\lambda}(Y)=\sum_{\mu\in W\lambda}Y^{\mu}$ ($\lambda\in
\Lambda^+$).

It follows from \eqref{lusztig} that
${\mathcal{A}}_Y^W$ lies in the center
${\mathcal{Z}}(H)$ of $H$. In fact, by \cite[prop. 3.11]{L} we know that
${\mathcal{Z}}(H)={\mathcal{A}}_Y^W$ (which also follows
from results in section 6).
The action of ${\mathcal{A}}_Y^W$ on
${\mathcal{A}}$ through
the Noumi representation $\pi_{\mathbf{t},q}$ preserves
${\mathcal{A}}^W$, compare e.g. with \cite[4.8]{M3}.

We extend the action of ${\mathcal{W}}$ on ${\mathcal{A}}$
to an action on the quotient field
${\mathbb{C}}(x)=\hbox{Quot}({\mathcal{A}})$
by requiring $w\in {\mathcal{W}}$ to be
an automorphism of ${\mathbb{C}}(x)$.
Let ${\mathbb{C}}(x)[{\mathcal{W}}] \subset
\hbox{End}_{\mathbb{C}}({\mathbb{C}}(x))$ be the subalgebra
generated by ${\mathbb{C}}(x)$ (acting as multiplication
operators) and by ${\mathcal{W}}$.
Observe that ${\mathcal{H}}\subset
{\mathbb{C}}(x)[{\mathcal{W}}]$, and that
\[{\mathbb{C}}(x)[{\mathcal{W}}]=\bigoplus_{w\in
{\mathcal{W}}}{\mathbb{C}}(x)w=\bigoplus_{w\in W, \lambda\in
\Lambda} {\mathbb{C}}(x)\tau(\lambda)w
\]
as a ${\mathbb{C}}(x)$-submodule of
$\hbox{End}_{\mathbb{C}}({\mathbb{C}}(x))$,
see the proof of \cite[thm. 3.2]{Sa}.
Furthermore, ${\mathbb{C}}(x)[\tau(\Lambda)]=\bigoplus_{\lambda\in \Lambda}
{\mathbb{C}}(x)\tau(\lambda)$ is the subalgebra of
${\mathbb{C}}(x)[{\mathcal{W}}]$
consisting of $q$-difference operators with coefficients in ${\mathbb{C}}(x)$.

With $D\in {\mathbb{C}}(x)[{\mathcal{W}}]$, say
\[ D=\sum_{w\in W}D(x,w)w,\qquad D(x,w)\in
{\mathbb{C}}(x)[\tau(\Lambda)],
\]
we associate a $q$-difference operator by
\[ D_{sym}=\sum_{w\in W}D(x,w)\in {\mathbb{C}}(x)[\tau(\Lambda)].
\]
Observe that $Df=D_{sym}f$ if $f\in {\mathbb{C}}(x)$ is $W$-invariant.
Proposition \ref{triangularY} and lemma \ref{compatibleactions}
imply that the $q$-difference operators $f(Y)_{sym}$ ($f\in {\mathcal{A}}^W$)
are triangular endomorphisms of ${\mathcal{A}}^W$:
\[
f(Y)_{sym}\,m_{\lambda}=f(\gamma_{\lambda})m_{\lambda}+
\sum_{\mu\in \Lambda^+: \mu<\lambda}c_{\lambda,\mu}m_{\mu},\qquad
f\in {\mathcal{A}}^W,\,\, \lambda\in \Lambda^+
\]
for certain constants $c_{\lambda,\mu}\in {\mathbb{C}}$. This
immediately implies the following result.
\begin{thm}\label{symmetric}
There exists a unique basis $\{P_{\lambda}^+\}_{\lambda\in \Lambda^+}$
of ${\mathcal{A}}^W$ such that
\begin{enumerate}
\item[{--}] $P_{\lambda}^+=m_{\lambda}+\sum_{\mu\in \Lambda^+: \mu<\lambda}
c_{\lambda,\mu}m_{\mu}$ for certain constants $c_{\lambda,\mu}$,
\item[{--}] $f(Y)_{sym}\,P_{\lambda}^+=f(\gamma_{\lambda})P_{\lambda}^+$ for
all $f(Y)\in {\mathcal{A}}_Y^W$,
\end{enumerate}
for all $\lambda\in \Lambda^+$.
\end{thm}
Noumi \cite{N}
identified the $q$-difference operator
\[ m_{\epsilon_1}(Y)_{sym}
=\bigl(Y_1+\cdots +Y_n+Y_1^{-1}+\cdots +Y_n^{-1}\bigr)_{sym}\in
{\mathbb{C}}(x)[\tau(\Lambda)]
\]
with Koornwinder's \cite{Ko} second order $q$-difference
operator.
Explicitly, Noumi \cite{N} showed that the
$q$-difference operator
$L=m_{\epsilon_1}(Y)_{sym}-m_{\epsilon_1}(\gamma_0)$ is given by
\begin{equation}\label{L}
L=\sum_{j=1}^n\bigl(\phi_j^+(x)(\tau(\epsilon_j)-1)+
\phi_j^-(x)(\tau(-\epsilon_j)-1)\bigr)
\end{equation}
with $\phi_j^-(x)=\phi_j^+(x_1^{-1},\ldots,x_n^{-1})$ and
\begin{equation*}
\begin{split}
\phi_j^{+}(x)=&(t_{0}t_{n})^{-1}t^{2(1-n)}
\frac{(1-ax_j)(1-bx_j)(1-cx_j)(1-dx_j)}{(1-x_j^2)(1-qx_j^2)}\\
&.\prod_{i\not=j}\frac{(1-t^2x_ix_j)(1-t^2x_i^{-1}x_j)}
{(1-x_ix_j)(1-x_i^{-1}x_j)}.
\end{split}
\end{equation*}
Here $\{a,b,c,d\}$ is related to the multiplicity
function $\mathbf{t}$ by
\begin{equation}\label{parametrization}
\{a,b,c,d\}=
\{t_{0}t_{0}^\vee q^{1/2},
-t_{0}(t_{0}^\vee)^{-1}q^{1/2},
t_{n}t_{n}^\vee,
-t_{n}(t_{n}^\vee)^{-1}\}.
\end{equation}
Since the spectrum of $L\in \hbox{End}_{\mathbf{C}}({\mathcal{A}}^W)$
is already simple (see \cite{Ko}),
this result implies that the $W$-invariant Laurent polynomials $P_{\lambda}^+$
($\lambda\in\Lambda^+$) coincide with Koornwinder's \cite{Ko}
multivariable analogues of the Askey-Wilson polynomials.
\begin{Def}
The $W$-invariant Laurent polynomial
$P_{\lambda}^+(\cdot)=P_{\lambda}^+(\cdot;\mathbf{t};q)$
\textup{(}$\lambda\in\Lambda^+$\textup{)} is called
the monic, symmetric Koornwinder polynomial of degree $\lambda\in
\Lambda^+$.
\end{Def}


\section{The action of the Hecke algebra of type $C_n$}

We associate a dual
multiplicity function $\tilde{\mathbf{t}}$ with the multiplicity
function $\mathbf{t}$ by interchanging the value of $\mathbf{t}$
on the ${\mathcal{W}}$-orbit ${\mathcal{W}}a_0$ with its value 
on the ${\mathcal{W}}$-orbit ${\mathcal{W}}a_n^\vee$. 
In other words, $\tilde{\mathbf{t}}$
is the unique multiplicity function of $S$ satisfying
\[
\tilde{t}_{0}=t_n^\vee, \qquad
\tilde{t}_0^\vee=t_0^\vee,\qquad
\tilde{t}=t,\qquad
\tilde{t}_{n}^\vee=t_{0},\qquad
\tilde{t}_{n}=t_{n}.
\]
We write $\tilde{\mathbf{k}}$ and $\tilde{\mathbf{k}}^\vee$ for the
associated multiplicity functions of $R$ and $R^\vee$ respectively, and
$\widetilde{v}_{\beta}(x)=v_{\beta}(x;\tilde{\mathbf{t}};q)$
for the function $v_{\beta}(\cdot)$ \eqref{v} with respect to dual
parameters.
Observe that Lusztig's formulas \eqref{lusztig} can now be
written in a uniform way:
\begin{equation}\label{lusztig2}
T_if(Y)-(s_if)(Y)T_i=\bigl((\tilde{t}_{a_i}-\tilde{t}_{a_i}^{-1})+
(\tilde{t}_{a_i/2}-\tilde{t}_{a_i/2}^{-1})Y^{-a_i/2}\bigr)
\left(\frac{f(Y)-(s_if)(Y)}{1-Y^{-a_i}}\right)
\end{equation}
for $i=1,\ldots,n$ and $f\in {\mathcal{A}}$.
In the following proposition we expand $T_iP_{\lambda}$
as a linear combination of non-symmetric Koornwinder polynomials.

\begin{prop}\label{symmnon}
Let $i\in \{1,\ldots,n\}$ and $\lambda\in \Lambda$. Then
\begin{equation}
\label{convex}
T_iP_{\lambda}=\xi_{i}(\gamma_{\lambda})P_{\lambda}+
\eta_{i}(\gamma_{\lambda})P_{s_i\lambda}
\end{equation}
with
\begin{equation}
\label{xi}
\xi_{i}(x)=\tilde{t}_i-\tilde{t}_i^{-1}\widetilde{v}_{-a_i}(x)
=\frac{\bigl(\tilde{t}_{a_i}^{-1}-
\tilde{t}_{a_i}\bigr)x^{a_i}+
\bigl(\tilde{t}_{a_i/2}^{-1}-\tilde{t}_{a_i/2}\bigr)x^{a_i/2}}{1-x^{a_i}}
\end{equation}
and
\begin{equation}\label{eta}
\eta_i(\gamma_{\lambda})=
\begin{cases}
\tilde{t}_i, &\hbox{ if } \langle \lambda,a_i\rangle<0,\\
\tilde{t}_{i}^{-3}\widetilde{v}_{a_i}(\gamma_{\lambda})
\widetilde{v}_{-a_i}(\gamma_{\lambda}),\qquad &\hbox{ if }
\langle \lambda,a_i\rangle\geq 0.
\end{cases}
\end{equation}
\end{prop}
\begin{proof}
The proof is based on the following consequence of
Lusztig's formula \eqref{lusztig2} and theorem
\ref{nonsymmetric}: let $\lambda\in \Lambda$ and $i\in \{1,\ldots,n\}$,
then
\begin{equation}\label{h}
\bigl(f(Y)-(s_if)(\gamma_{\lambda}\bigr)\bigr)T_iP_{\lambda}=
\bigl(f(\gamma_{\lambda})-(s_if)(\gamma_{\lambda}\bigr)\bigr)
\xi_i(\gamma_{\lambda})P_{\lambda}, \qquad \forall f\in
{\mathcal{A}},
\end{equation}
with $\xi_i$ given by \eqref{xi}.

Suppose now first that $\langle \lambda, a_i\rangle=0$, i.e. that
$s_i\lambda=\lambda$. Using remark \ref{casestable}, we see that
$\xi_i(\gamma_{\lambda})+\eta_i(\gamma_{\lambda})=t_i$, so we have
to show that $T_iP_{\lambda}=t_iP_{\lambda}$. Now \eqref{h},
theorem \ref{nonsymmetric} and remark \ref{casestable} imply
that $T_iP_{\lambda}$ is a constant multiple of $P_{\lambda}$.
The constant multiple can be determined by computing the leading coefficient of
$T_iP_{\lambda}$ using the identity
$T_i=s_i{\mathcal{R}}(a_i)^{-1}+t_i-t_i^{-1}$ and using
lemma \ref{triangularR}.

Suppose now that $\langle \lambda,a_i\rangle\not=0$, then
\eqref{h}, lemma \ref{compatibleactions} and theorem
\ref{nonsymmetric} imply that $T_iP_{\lambda}$ is of the form
\eqref{convex} for some constant $\eta_i(\gamma_{\lambda})$,
with $\xi_i$ given by \eqref{xi}.
If $\langle \lambda,a_i\rangle<0$, then leading term
considerations using lemma \ref{triangularR} and
theorem \ref{nonsymmetric} show that
$\eta_i(\gamma_{\lambda})=t_i=\tilde{t}_i$.
The expression for $\eta_i(\gamma_{\lambda})$
when $\langle \lambda,a_i\rangle>0$ follows now easily
by applying $T_i$ on both sides of \eqref{convex} and using
the quadratic relation $(T_i-t_i)(T_i+t_i^{-1})=0$,
compare with \cite[prop. 4.1]{NS} for the proof in the rank one setting.
\end{proof}

We write $S_i=[T_i,Y^{a_i}]=T_iY^{a_i}-Y^{a_i}T_i\in {\mathcal{H}}$
($i=1,\ldots,n$). It follows from \eqref{lusztig2} that the
$S_i$ satisfy the fundamental commutation relations
$f(Y)S_i=S_i(s_if)(Y)$ for all $f\in {\mathcal{A}}$, cf.
\cite[sect. 5]{Sa}.
We call $S_i$ the (non-affine) intertwiner associated with
the simple reflection $s_i$.

We use here a slightly different definition for the intertwiners $S_i$
compared with Sahi's \cite{Sa}, \cite{Sa2} intertwiners.
The advantage of the present definition is that
the $S_i$ ($i=1,\ldots,n$) satisfy the $C_n$-braid
relations, see remark \ref{braidS}.
In particular, we may write $S_w=S_{i_1}\cdots S_{i_r}$ for a
reduced expression $w=s_{i_1}\cdots s_{i_r}\in W$, and $S_w$
satisfies the intertwining property $S_wf(Y)=(wf)(Y)S_w$ for
all $f\in {\mathcal{A}}$. See also the paper \cite{NUKW}, 
in which yet another definition for the  intertwiners $S_i$ ($i=0,\ldots,n$)
is used (including a non-affine intertwiner $S_0$). The intertwiners
in \cite{NUKW}, which 
satisfy the $\widetilde{C}_n$-braid relations, are
used to prove a Rodrigues type formula for
non-symmetric Koornwinder polynomials.

With our present  conventions,  the action of the intertwiners 
on the non-sym\-me\-tric Koornwinder
polynomials is easily determined from proposition \ref{symmnon}.
We give here only the action of the intertwiner $S_i$
corresponding to a simple reflection $s_i$ of $W$.

\begin{cor}\label{normintertwiner}
$S_iP_{\lambda}=\bigl(\gamma_{\lambda}^{a_i}-
\gamma_{s_i\lambda}^{a_i}\bigr)
\eta_i(\gamma_{\lambda})P_{s_i\lambda}$
for $i=1,\ldots,n$ and $\lambda\in \Lambda$.
\end{cor}

\begin{rem}
Proposition \ref{symmnon} refines the 
non-affine part of Sahi's recursion formula
\cite[thm. 18]{Sa2}, while corollary \ref{normintertwiner} can be seen as
a refinement of the non-affine part of \cite[thm. 5.3]{Sa}.
Indeed, in \cite[thm. 18]{Sa2} and \cite[thm. 5.3]{Sa} the
formulas are given up to an unknown multiple constant.
The unknown constants in the affine part of \cite[thm. 18]{Sa2} and
\cite[thm. 5.3]{Sa} can also be computed, but this requires 
the duality properties and the evaluation formulas for 
the non-symmetric Koornwinder polynomials, see proposition
\ref{transfer} and theorem \ref{evaluationthm}.
\end{rem}

By theorem \ref{nonsymmetric}, proposition \ref{symmnon} and corollary
\ref{normintertwiner} we have a complete description of the
action of $H(R;\mathbf{k})$ on the non-symmetric Koornwinder polynomials
under the Noumi representation $\pi_{\mathbf{t},q}$. {}From this
the $H$-module structure of ${\mathcal{A}}$ can be described in
detail, see also Sahi \cite{Sa}. The result is as follows. We write
\begin{equation}
\label{AH}
{\mathcal{A}}=\bigoplus_{\lambda\in \Lambda^+}
{\mathcal{A}}(\lambda),\qquad {\mathcal{A}}(\lambda)=
\hbox{span}\{P_{\mu} \, | \, \mu\in W\lambda \}.
\end{equation}
Recall that the parameters $\mathbf{t}$ and $q$ are assumed to be generic.
\begin{thm}\label{noumithm}
The direct sum decomposition \eqref{AH}
is the multiplicity-free, irreducible
decomposition of ${\mathcal{A}}$ as a
$(\pi_{\mathbf{t},q}, H(R;\mathbf{k}))$-module.
Furthermore, \eqref{AH} is
the decomposition of ${\mathcal{A}}$ into
isotypical components under the action of the center
${\mathcal{Z}}(H(R;\mathbf{k}))={\mathcal{A}}_Y^W$.
The central character $\chi_{\lambda}$
of ${\mathcal{A}}(\lambda)$ is given by
$\chi_{\lambda}(f)=f\bigl(\gamma_{\lambda})$ for
$f\in {\mathcal{Z}}(H(R;\mathbf{k}))$.
\end{thm}

We end this section by defining anti-symmetric Koornwinder polynomials and by
expanding (anti-)symmetric Koornwinder
polynomials in terms of non-symmetric Koornwinder polynomials.

We associate 
a function $\{t_w\}_{w\in {\mathcal{W}}}$ 
with the multiplicity function $\mathbf{t}=\{t_{\beta}\}_{\beta\in S}$
by defining $t_w=t_{i_1}\ldots t_{{i_r}}$ for a reduced expression
$w=s_{i_1}\cdots s_{i_r}\in {\mathcal{W}}$.
\begin{rem}\label{compatible}
Recall
the well-known fact that
$R^+\cap w^{-1}R^-=\{ \beta_1,\ldots, \beta_r\}$
with the $r$ distinct positive roots $\beta_j$ given by
\[ \beta_j=s_{i_r}\cdots s_{i_{j+1}}a_{i_j}\quad
(j=1,\ldots,r-1),\,\,\, \beta_r=a_{i_r}.
\]
In particular, it follows that
\[
t_w=\prod_{\beta\in R^+\cap w^{-1}R^-}t_{\beta},\qquad w\in {\mathcal{W}}.
\]
Restricted to the finite Weyl group $W$, the expression for $t_w$
reduces to
\[t_w=\prod_{\alpha\in \Sigma^+\cap w^{-1}\Sigma^-}t_{\alpha},\qquad
w\in W.
\]
Observe in particular that $\tilde{t}_w=t_w$ when $w\in W$, where
$\{\tilde{t}_{w}\}_{w\in {\mathcal{W}}}$ is the function 
associated with the dual multiplicity function $\tilde{\mathbf{t}}$.
\end{rem}
Let $\chi_{\pm}: H_0\rightarrow {\mathbb{C}}$ be the trivial and
alternating character of the Hecke algebra $H_0$, i.e.
$\chi_{\pm}(T_i)=\pm t_{i}^{\pm 1}$ for $i=1,\ldots,n$.
Then the corresponding mutually orthogonal, primitive idempotents
are given by
\begin{equation}\label{C}
C_\pm=\frac{1}{\sum_{w\in W}t_w^{\pm 2}}\sum_{w\in W}(\pm 1)^{l(w)}t_w^{\pm
1}T_w.
\end{equation}
We define for $\lambda\in \Lambda^+$,
\begin{equation}\label{apm}
\begin{split}
{\mathcal{A}}_{\pm}(\lambda)&=\{ f\in {\mathcal{A}}(\lambda) \,
| \, C_{\pm}f=f \}\\
&=\{ f\in {\mathcal{A}}(\lambda) \, | \, (T_i\mp t_i^{\pm
1})f=0\,\,\, \forall i=1,\ldots,n\}.
\end{split}
\end{equation}
Observe in particular that
${\mathcal{A}}_+(\lambda)={\mathcal{A}}(\lambda)\cap
{\mathcal{A}}^W$.
Let $\Lambda^{++}=\kappa+\Lambda^+$ be the cone of regular dominant weights, 
where
\begin{equation}\label{kappa}
\kappa=\frac{1}{2}\sum_{\alpha\in\Sigma^+}\alpha=\sum_{i=1}^n\omega_i.
\end{equation}
\begin{thm}\label{antisymmetric}
{\bf (i)} ${\mathcal{A}}_+(\lambda)$ is spanned by the symmetric Koornwinder
polynomial $P_{\lambda}^+$ for all $\lambda\in \Lambda^+$.

{\bf (ii)}
${\mathcal{A}}_-(\lambda)$ is one-dimensional if $\lambda\in \Lambda^{++}$
and zero dimensional otherwise.
For $\lambda\in \Lambda^{++}$ there exists
a unique $P_{\lambda}^{-}\in {\mathcal{A}}_{-}(\lambda)$
of the form $P_{\lambda}^{-}=x^{\lambda}+\sum_{\mu\prec
\lambda}e_{\mu}x^{\mu}$ for certain constants $e_{\mu}\in {\mathbb{C}}$.

{\bf (iii)}
The expressions for $P_{\lambda}^{+}$
\textup{(}$\lambda\in\Lambda^+$\textup{)} and for $P_{\lambda}^-$
\textup{(}$\lambda\in\Lambda^{++}$\textup{)} as linear combinations
of the non-sym\-me\-tric Koornwinder polynomials $P_{\mu}$ \textup{(}$\mu\in
W\lambda$\textup{)} are given by
\[ P_{\lambda}^{\pm}=\sum_{\mu\in
W\lambda}c_{\lambda,\mu}^{\pm}P_{\mu}
\]
with the coefficients $c_{\lambda,\mu}^{\pm}$ 
\textup{(}$\mu\in W\lambda$\textup{)} given by
\begin{equation}\label{clm}
\begin{split}
c_{\lambda,\mu}^{\pm}&=(\pm 1)^{l(w_{\mu})}\tilde{t}_{w_{\mu}}^{-2}
\prod_{\stackrel{\alpha\in\Sigma^+}{\langle
\mu,\alpha\rangle<0}}\widetilde{v}_{\pm\alpha}(\gamma_{\mu})\\
&=(\pm 1)^{l(w_{\mu})}
\tilde{t}_{w_{\mu}}^{-2}\prod_{\alpha\in \Sigma^+\cap
w_{\mu}^{-1}\Sigma^-}\widetilde{v}_{\pm \alpha}\bigl(\gamma_{\lambda}^{-1}),
\end{split}
\end{equation}
where $w_{\mu}$ is the element of minimal length in $W$ such that
$w_{\mu}\lambda(=w_{\mu}\mu^+)=\mu$.
\end{thm}
\begin{proof}
We first introduce some notations and deduce some preliminary results
which we will need for the proof.

Let $\lambda\in \Lambda^+$ and let
$Q_{\lambda}^{\pm}=\sum_{\mu\in
W\lambda}d_{\lambda,\mu}^{\pm}P_{\mu}$
be an element in ${\mathcal{A}}_{\pm}(\lambda)$.
By proposition \ref{symmnon} we 
have for $\mu\in \Lambda$ and $i\in\{1,\ldots,n\}$
that
\[\bigl(T_i\mp t_{i}^{\pm 1}\bigr)P_{\mu}=
\xi_{i}^{\pm}(\gamma_{\mu})P_{\mu}+
\eta_{i}(\gamma_{\mu})P_{s_i\mu}
\]
with
\begin{equation}\label{tussen}
\xi_{i}^{\pm}(x)=\xi_i(x)\mp t_i^{\pm 1}=
\mp \tilde{t}_{i}^{-1}\widetilde{v}_{\mp a_i}(x).
\end{equation}
It follows now from the relations $(T_i\mp t_{i}^{\pm 1})Q_{\lambda}^{\pm}=0$
that the coefficients $d_{\lambda,\mu}^{\pm}$ satisfy the
recurrence relations
\begin{equation}
\label{recurrenc}
d_{\lambda,s_i\mu}^{\pm}=
\left(-\frac{\xi_{i}^{\pm}(\gamma_{\mu})}{\eta_{i}(\gamma_{s_i\mu})}
\right)d_{\lambda,\mu}^{\pm}
\end{equation}
for $\mu\in W\lambda$ and $i\in \{1,\ldots,n\}$ such that
$\langle \mu,a_i\rangle\not=0$.
Let now $\mu\in W\lambda$ and choose a reduced expression $w_{\mu}=s_{i_1}\cdots
s_{i_r}$. We set
\[ \mu_j=s_{i_{j+1}}\cdots s_{i_r}\lambda \quad
(j=0,\ldots,r-1),\quad \mu_r=\lambda.
\]
Then $\langle \mu_j,a_{i_j}\rangle>0$ for $j=1,\ldots,r$ by remark
\ref{compatible}.
Iterating \eqref{recurrenc}, we thus obtain
\begin{equation}
\label{coeffnew}
d_{\lambda,\mu}^{\pm}=d_{\lambda,\lambda}^{\pm}
\prod_{j=1}^r\left(-\frac{\xi_{i_j}^{\pm}(\gamma_{\mu_j})}
{\eta_{i_j}(\gamma_{\mu_{j-1}})}\right).
\end{equation}

{\it Proof of} {\bf (i)} The recurrence formula \eqref{coeffnew}
implies that $\hbox{dim}({\mathcal{A}}_{+}(\lambda))\leq 1$
for all $\lambda\in \Lambda^+$.
On the other hand, theorem \ref{symmetric} implies that
the symmetric Koornwinder polynomial $P_{\lambda}^+$ is a non-zero
element in ${\mathcal{A}}_+(\lambda)$ for all $\lambda\in \Lambda^+$.

{\it Proof of }{\bf (ii)} Again by \eqref{coeffnew}, we have
$\hbox{dim}({\mathcal{A}}_{-}(\lambda))\leq 1$.
Let $\lambda\in \Lambda^+\setminus \Lambda^{++}$.
Let $s_i\in W$ be a simple reflection in $W$
which stabilizes $\lambda$. Then proposition
\ref{symmnon} implies that $T_iP_{\lambda}=t_{i}P_{\lambda}$.
Hence the coefficient of $P_{\lambda}$
in the expansion of $(T_i+t_{i}^{-1})Q_{\lambda}^-$ as a linear
combination of the $P_{\nu}$ ($\nu\in W\lambda$), is
$(t_{i}+t_{i}^{-1})d_{\lambda,\lambda}^-$. On the other hand,
$(T_i+t_{i}^{-1})Q_{\lambda}^-=0$, hence we conclude that
$d_{\lambda,\lambda}^-=0$. Then \eqref{coeffnew} implies $d_{\lambda,\mu}^-=0$
for all $\mu\in W\lambda$, hence $Q_{\lambda}^-=0$. This proves
that ${\mathcal{A}}_-(\lambda)=\{0\}$ if $\lambda\in \Lambda^+\setminus
\Lambda^{++}$.

Let now $\lambda\in \Lambda^{++}$. Then $C_-P_{\sigma\lambda}\in
{\mathcal{A}}_-(\lambda)$, where $\sigma\in W$ is the 
longest Weyl group element.
Furthermore, by proposition
\ref{symmnon} we have 
$C_-P_{\sigma\lambda}=\sum_{\mu\in W\lambda}d_{\mu}P_{\mu}$
with $d_{\lambda}=(\sum_{w\in
W}t_w^{-2})^{-1}(-1)^{l(\sigma)}\not=0$,
so $C_-P_{\sigma\lambda}$ is a non-zero 
element in ${\mathcal{A}}_-(\lambda)$. Hence
$\hbox{dim}({\mathcal{A}}_-(\lambda))=1$ if $\lambda\in \Lambda^{++}$.
The triangularity statement follows from theorem \ref{nonsymmetric} 
and from the
fact that the coefficient $d_{\lambda}$ in the above expansion of
$C_-P_{\sigma\lambda}$ is non-zero.

{\it Proof of }{\bf (iii)}
We have to show that
\begin{equation}
\label{remain}
d_{\lambda,\mu}^{\pm}=d_{\lambda,\lambda}^{\pm}c_{\lambda,\mu}^{\pm},\qquad
\mu\in W\lambda,
\end{equation}
where the $d_{\lambda,\mu}^{\pm}$ are the expansion coefficients of
$Q_{\lambda}^{\pm}$ in terms of non-symmetric Koornwinder polynomials $P_{\mu}$
($\mu\in W\lambda$), and where $c_{\lambda,\mu}^{\pm}$ is given by
\eqref{clm}. We use again the recurrence formula \eqref{coeffnew}
for the coefficients $d_{\lambda,\mu}^{\pm}$.
We set
\[ \alpha_1=-a_{i_1},\quad
\alpha_j=-s_{i_1}\cdots s_{i_{j-1}}a_{i_j} \quad
(j=2,\ldots,r),
\]
then the $\alpha_j$ ($j=1,\ldots,r$) are mutually different
and
\begin{equation}\label{char}
\{\alpha_1,\ldots,\alpha_r\}=\Sigma^-\cap w_{\mu}\Sigma^+
=\{ \alpha\in \Sigma^- \, | \, \langle
\mu,\alpha\rangle>0 \},
\end{equation}
see remark \ref{compatible}. Now observe that
$\langle \mu_{j-1},a_{i_{j}}\rangle=-\langle \mu_{j},
a_{i_{j}}\rangle<0$ for all $j=1,\ldots,r$, so that
\begin{equation}\label{beta}
\eta_{i_{j}}(\gamma_{\mu_{j-1}})=
\tilde{t}_{\alpha_{j}}
\qquad (j=1,\ldots,r).
\end{equation}
Furthermore, observe that
\begin{equation}
\label{innok}
(\gamma_{\mu_j})^{a_{i_j}}=(\gamma_{\mu})^{\alpha_j}=
(\gamma_{\lambda})^{w_{\mu}^{-1}\alpha_j},\qquad
j=1,\ldots,r
\end{equation}
by lemma \ref{compatibleactions}.
Substituting \eqref{beta} and \eqref{tussen}
in \eqref{coeffnew} and using \eqref{innok} and the characterization of the
roots $\{\alpha_j\}_{j=1}^r$ (see \eqref{char}), we obtain
\eqref{remain}.
\end{proof}

\begin{Def}
The Laurent polynomial $P_{\lambda}^-\in {\mathcal{A}}_-(\lambda)$
\textup{(}$\lambda\in \Lambda^{++}$\textup{)} is
called the anti-symmetric Koornwinder polynomial of degree $\lambda$.
\end{Def}

\begin{rem}
{\bf (i)} Part {\bf (i)} of theorem \ref{antisymmetric} was also observed by
Sahi \cite[cor. 6.6]{Sa}.

{\bf (ii)} Theorem \ref{antisymmetric} extends Macdonald's 
\cite[sect. 6]{M3} explicit expansion formulas for
the (anti-)symmetric Macdonald
polynomials associated with root systems 
of classical type.
\end{rem}


\section{Spectral difference-reflection operators and duality}

Let 
$x_{\lambda}=\gamma_{\lambda}(\tilde{\mathbf{k}},q)$ ($\lambda\in \Lambda$)
be the spectrum \eqref{gl} of the $\widetilde{Y}$-operators, and denote
$\widetilde{\mathcal{H}}={\mathcal{H}}(S;\tilde{\mathbf{t}};q)$ 
for the double affine Hecke algebra with respect to dual parameters.
We define
evaluation mappings $\hbox{Ev}: {\mathcal{H}}
\rightarrow {\mathbb{C}}$ and $\widetilde{\hbox{Ev}}: \widetilde{\mathcal{H}}
\rightarrow {\mathbb{C}}$ by
\[ \hbox{Ev}(X)=\bigl(X(1)\bigr)(x_0^{-1}),\qquad
\widetilde{\hbox{Ev}}(\widetilde{X})=
\bigl(\widetilde{X}(1)\bigr)
(\gamma_0^{-1})
\]
for $X\in {\mathcal{H}}$ and $\widetilde{X}\in
\widetilde{\mathcal{H}}$,
where $1\in {\mathcal{A}}$ is the Laurent polynomial identically equal to one.

The evaluation $\hbox{Ev}(P_{\lambda}(z))=P_{\lambda}(x_0^{-1})$
of the non-symmetric monic Koornwinder polynomial $P_{\lambda}(\cdot)=
P_{\lambda}(\cdot;\mathbf{t};q)$
is generically non-zero by
the analytic dependence of $P_{\lambda}$ on $\mathbf{t}$ and $q$
(we use here that $P_{\lambda}(x)=x^{\lambda}$
when $t_a=1$ for all $a\in S$). Similarly, $\hbox{Ev}(P_{\lambda}^+(z))=
P_{\lambda}^+(x_0^{\pm 1})$ is non-zero for generic parameter values $\mathbf{t}$
and $q$. In section 9 we explicitly evaluate
$P_{\lambda}(x_0^{-1})$ and $P_{\lambda}^+(x_0)=P_{\lambda}^+(x_0^{-1})$,
so that the generic conditions on the parameters can be made
completely explicit.

\begin{Def}
{\bf (i)} Let $E({\gamma_{\lambda}};\cdot)=
E(\gamma_{\lambda};\cdot;\mathbf{t};q)$ be the constant multiple
of the non-symmetric Koornwinder polynomial $P_{\lambda}(\cdot)$
of degree $\lambda\in \Lambda$ which takes the value one at $x=x_0^{-1}$.\\
{\bf (ii)} Let
$E^+(\gamma_{\lambda};\cdot)=E^+({\gamma_{\lambda}};\cdot;\mathbf{t};q)$
be the constant multiple of the
symmetric Koornwinder polynomial $P_{\lambda}^+(\cdot)$
of degree $\lambda\in \Lambda^+$ which takes the value one at $x=x_0$.
\end{Def}

Sahi \cite{Sa} showed that the role of the geometric parameter
$x=x_{\mu}$ and of the spectral parameter $\gamma=\gamma_{\lambda}$
are (in a suitable sense) interchangeable for the
renormalized Koornwinder polynomials
$E(\gamma;x^{-1})$ and $E^+(\gamma;x)$, see also van Diejen \cite{vD2}
for a sub-class of the symmetric Koornwinder polynomials.
These duality properties stem from a particular anti-algebra isomorphism
of the double affine Hecke algebra ${\mathcal{H}}$, which we define
now first.

Recall the notations 
$T_0^\vee=T_0^{-1}z^{-a_0^\vee}\in {\mathcal{H}}$
and $T_n^\vee=z^{-a_n^\vee}T_n^{-1}\in {\mathcal{H}}$
for the simple generators associated with $a_0^\vee$ and $a_n^\vee$
respectively, see remark \ref{Spresent}. We set
\[
U_n= T_1T_2\cdots T_{n-1}T_n^\vee T_{n-1}^{-1}\cdots T_2^{-1}T_1^{-1},
\]
which is a conjugate of $T_n^\vee$ in ${\mathcal{H}}(S;\mathbf{t};q)$.

Set $\mathbf{t}^{-1}=(t_\beta^{-1})_{\beta\in S}$ for the inverse of
the multiplicity function $\mathbf{t}$. We write
$T_i^\prime, T_j^{\vee\,\prime}$, $U_n^\prime$, $Y^{\prime\,\lambda}$
$z^{\prime\,\lambda}$ for the elements $T_i,T_j^\vee$, $U_n$, $Y^{\lambda}$
and $z^{\lambda}$ in the double affine Hecke algebra
${\mathcal{H}}^\prime={\mathcal{H}}(S;\mathbf{t}^{-1};q^{-1})$.
Similarly, we write $\widetilde{T}_i,\ldots$ (respectively
$\widetilde{T}_i^\prime,\ldots$) for the elements $T_i,\ldots$
in the double affine Hecke algebra $\widetilde{{\mathcal{H}}}$
(respectively ${\widetilde{\mathcal{H}}}^\prime=
{\mathcal{H}}(S;{\tilde{\mathbf{t}}}^{-1};q^{-1})$).
The following theorem was proved by Sahi \cite[thm. 4.2]{Sa}.
\begin{thm}\label{epsilon}
There exists a unique algebra isomorphism
$\epsilon=\epsilon_{\mathbf{t},q}: {\mathcal{H}}\rightarrow
{\widetilde{\mathcal{H}}}^\prime$ satisfying
$\epsilon(T_0)=(\widetilde{U}_n^\prime)^{-1}$,
$\epsilon(z_i)=\widetilde{Y}_i^\prime$
and $\epsilon(T_i)=(\widetilde{T}_i^\prime)^{-1}$ for $i=1,\ldots,n$.
Furthermore, 
$\epsilon_{\mathbf{t},q}^{-1}=\epsilon_{\tilde{\mathbf{t}}^{-1},q^{-1}}$.
\end{thm}

The isomorphism $\epsilon$ is a crucial building block for
Sahi's \cite{Sa} duality anti-iso\-mor\-phism of the double affine Hecke
algebra ${\mathcal{H}}$. In fact, the duality 
anti-isomorphism is obtained by
composing $\epsilon$ with the anti-isomorphism $\ddagger$
defined in the following lemma.

\begin{lem}\label{adj}
There exists a unique algebra isomorphism $\dagger=\dagger_{\mathbf{t},q}:
{\mathcal{H}}\rightarrow {\mathcal{H}}^\prime$
\textup{(}respectively anti-algebra isomorphism
$\ddagger=\ddagger_{\mathbf{t},q}:
{\mathcal{H}}\rightarrow {\mathcal{H}}^\prime$\textup{)}
satisfying $T_i\mapsto (T_i^\prime)^{-1}$ \textup{(}$i=0,\ldots,n$\textup{)}
and $z_j\mapsto (z_j^\prime)^{-1}$ \textup{(}$j=1,\ldots,n$\textup{)}.
\end{lem}
\begin{proof}
This follows directly from the presentation of ${\mathcal{H}}$
as given by Sahi \cite[sect. 3]{Sa}.
\end{proof}
\begin{rem}
{\bf (i)} Lemma \ref{adj} for $\ddagger$ was observed by Sahi
\cite[prop. 7.1]{Sa}.\\
{\bf (ii)} In proposition \ref{adjointprop} we interpret the anti-algebra
isomorphism $\ddagger$
as a $\ast$-structure on ${\mathcal{H}}\subset
\hbox{End}_{\mathbb{C}}({\mathcal{A}})$ induced from a
suitable non-degenerate bilinear form on ${\mathcal{A}}$.
\end{rem}

We write $\widetilde{\dagger}^\prime$ (respectively
$\widetilde{\ddagger}^\prime$)
for $\dagger$ (respectively $\ddagger$) with respect to the
parameters $(\tilde{\mathbf{t}}^{-1},q^{-1})$.

\begin{Def}
{\bf (i)} The algebra isomorphism
$\Phi=\Phi_{\mathbf{t},q}=\widetilde{\dagger}^\prime\circ\epsilon: {\mathcal{H}}
\rightarrow \widetilde{\mathcal{H}}$ is called the duality isomorphism of
${\mathcal{H}}$.\\
{\bf (ii)} The anti-algebra isomorphism
$\Psi=
\Psi_{\mathbf{t},q}=\widetilde{\ddagger}^\prime\circ\epsilon: {\mathcal{H}}
\rightarrow \widetilde{\mathcal{H}}$ is called
the duality anti-isomorphism of ${\mathcal{H}}$.
\end{Def}

Observe that $\Phi$ (respectively $\Psi$)  is uniquely characterized as the
(anti-)algebra homomorphism ${\mathcal{H}}\rightarrow \widetilde{\mathcal{H}}$
which maps $U_n$ to
$\widetilde{T}_0$, $T_i$ to
$\widetilde{T}_i$ and
$Y_i$ to $\widetilde{z}_i^{-1}$ for $i=1,\ldots,n$.
By \cite[sect. 7]{Sa}, the inverse of
$\Psi=\Psi_{\mathbf{t},q}$ is given by
$\widetilde{\Psi}=\Psi_{\tilde{\mathbf{t}},q}$.
\begin{rem}\label{braidS}
Observe that the image of the non-affine intertwiners
$S_i=[T_i, Y^{a_i}]\in {\mathcal{H}}$ ($i=1,\ldots,n$) under $\Psi$
is given by
\[
\Psi(S_i)=\tilde{t}_i^{-1}(\widetilde{z}^{-a_i}-\widetilde{z}^{a_i})
\widetilde{v}_{a_i}(\widetilde{z})s_i\in
\widetilde{\mathcal{H}}
\]
in view of the explicit expression for $\widetilde{T}_i$ (see 
theorem \ref{noumi}).
In particular, $\Psi(S_i)$ is of the
form $f_i(\widetilde{z}^{a_i^\vee})s_i$, with $f_i$ a Laurent
polynomial in one variable and with $f_1=f_2=\cdots =f_{n-1}$. {}From these facts
it is easy to prove
that $(S_1,\ldots,S_n)$ satisfies the $C_n$-braid relations
in ${\mathcal{H}}$.
\end{rem}

The two evaluation mappings $\hbox{Ev}$
and $\widetilde{\hbox{Ev}}$ are related via the duality anti-iso\-mor\-phism:
\[ \widetilde{\hbox{Ev}}\bigl(\Psi(X)\bigr)=\hbox{Ev}\bigl(X\bigr),\qquad
X\in {\mathcal{H}},
\]
see \cite[thm. 7.3]{Sa}.
This implies that the two pairings 
$B: {\mathcal{H}}\times \widetilde{\mathcal{H}}
\rightarrow {\mathbb{C}}$ and $\widetilde{B}: \widetilde{\mathcal{H}}\times
{\mathcal{H}}\rightarrow {\mathbb{C}}$ defined by
$B(X,\widetilde{X})=\hbox{Ev}\bigl(\widetilde{\Psi}(\widetilde{X})X\bigr)$
and $\widetilde{B}(\widetilde{X},X)=
\widetilde{\hbox{Ev}}\bigl(\Psi(X)\widetilde{X}\bigr)$ for $X\in {\mathcal{H}}$
and $\widetilde{X}\in \widetilde{\mathcal{H}}$ satisfy the duality
property
\begin{equation}
\label{dualityB}
B(X,\widetilde{X})=\widetilde{B}(\widetilde{X},X),\qquad X\in
{\mathcal{H}},\,\, \widetilde{X}\in \widetilde{\mathcal{H}}.
\end{equation}
Before we recall how \eqref{dualityB} implies
the duality properties of the Koornwinder polynomials,
we first collect some elementary identities for the
bilinear form $B$. The proof of the lemma
is similar to the proof in the
rank one setting, see \cite[lem. 10.5]{NS}.
\begin{lem}\label{easy}
Let $f\in {\mathcal{A}}$.
Let $X,X_1,X_2\in {\mathcal{H}}$ and
$\widetilde{X}, \widetilde{X}_1, \widetilde{X}_2\in
\widetilde{\mathcal{H}}$.

{\bf (i)} $B\bigl(X_1X_2,\widetilde{X}\bigr)=
B\bigl(X_2,\Psi(X_1)\widetilde{X}\bigr)$ and
$B(X, \widetilde{X}_1\widetilde{X}_2)=B(\widetilde{\Psi}(\widetilde{X}_1)X,
\widetilde{X}_2)$.

{\bf (ii)}
$B\bigl(XT_i,\widetilde{X}\bigr)=t_{i}B\bigl(X,\widetilde{X})$
for $i=0,\ldots,n$.

{\bf (iii)}
$B\bigl((X(f))(z),\widetilde{X}\bigr)=B(Xf(z),\widetilde{X})$
and $B\bigl(X,(\widetilde{X}(f))(\widetilde{z})\bigr)=
B\bigl(X,\widetilde{X}f(\widetilde{z})\bigr)$, where $(X(f))(z)$
is the multiplication operator in ${\mathcal{H}}$ corresponding to the
Laurent polynomial $X(f)\in {\mathcal{A}}$, and 
$Xf(z)$ is the product of the elements $X$ and 
$f(z)$ in ${\mathcal{H}}$.

\end{lem}

We write $\widetilde{E}({x_{\lambda}};\cdot)$ for the renormalized
non-symmetric Koornwinder polynomial
$E({x_{\lambda}};\cdot;\tilde{\mathbf{t}};q)$,
and similarly for $\widetilde{E}^+({x_{\lambda}};\cdot)$.
Observe now that by lemma \ref{easy} and theorem
\ref{nonsymmetric},
\begin{equation}\label{fgexp}
f(\gamma_{\lambda}^{-1})=\widetilde{B}\bigl(f(\widetilde{z}),
E({\gamma_{\lambda}};z)\bigr),\qquad g(x_{\mu}^{-1})=
B\bigl(g(z), \widetilde{E}({x_{\mu}};\widetilde{z})\bigr)
\end{equation}
for $f,g\in {\mathcal{A}}$ and $\lambda,\mu\in \Lambda$.
Taking $f=\widetilde{E}({x_{\mu}};\cdot)$ and $g=E({\gamma_{\lambda}};\cdot)$
and using the duality \eqref{dualityB} for the pairing, we arrive
at
\begin{equation}\label{duality}
E({\gamma_{\lambda}};x_{\mu}^{-1})=
\widetilde{E}({x_{\mu}};\gamma_{\lambda}^{-1}),
\qquad \lambda,\mu\in \Lambda
\end{equation}
which is the duality for the renormalized Koornwinder polynomials, see
\cite[thm. 7.4]{Sa}.
Similarly, we derive from theorem \ref{symmetric} and \eqref{dualityB} that
\begin{equation}
\label{dualitysymm}
E^+({\gamma_{\lambda}};x_{\mu})=\widetilde{E}^+({x_{\mu}};\gamma_{\lambda}),
\qquad \lambda,\mu\in \Lambda^+,
\end{equation}
see \cite[cor. 7.5]{Sa}.
Using the duality \eqref{duality}, we can rewrite
the action of $T_i$ ($i=1,\ldots,n$) and $U_n$ on the renormalized
Koornwinder polynomials $E({\gamma};\cdot)$
in terms of difference-reflection operators acting on the
spectral parameter $\gamma\in \hbox{Spec}(Y)=
\{\gamma_{\lambda} \, | \, \lambda\in \Lambda\}$. 
Define an action of ${\mathcal{W}}$ on $\hbox{Spec}(Y)$ by
$w\gamma_{\lambda}=\gamma_{w.\lambda}$ ($\lambda\in \Lambda$, $w\in
{\mathcal{W}}$).

\begin{prop}\label{transfer}
{\bf (i)} For $\gamma\in \hbox{Spec}(Y)$ we have
\[
\bigl(U_nE({\gamma};\cdot)\bigr)(x)=
\tilde{t}_{0}E({\gamma};x)+\tilde{t}_{0}^{-1}
\widetilde{v}_{a_0}(\gamma^{-1})\bigl(E(s_0\gamma;x)-E(\gamma;x)\bigr).
\]
{\bf (ii)} For $i=1,\ldots,n$ and $\gamma\in \hbox{Spec}(Y)$ we have
\[
\bigl(T_iE(\gamma;\cdot)\bigr)(x)=\tilde{t}_{i}E({\gamma};x)
+\tilde{t}_{i}^{-1}\widetilde{v}_{a_i}(\gamma^{-1})
\bigl(E({s_i\gamma};x)-E({\gamma};x)\bigr).
\]

\end{prop}
\begin{proof}
By \eqref{fgexp} and lemma \ref{easy} we have
\begin{equation*}
B\bigl(E(\gamma_{\lambda};z),
\widetilde{T}_i\widetilde{E}(x_{\mu};\widetilde{z})\bigr)=
\begin{cases}
\bigl(U_nE(\gamma_{\lambda};\cdot)\bigr)(x_{\mu}^{-1})\qquad
&\hbox{ if } i=0\\
\bigl(T_iE(\gamma_{\lambda};\cdot)\bigr)(x_{\mu}^{-1})\qquad
&\hbox{ if } i=1,\ldots,n
\end{cases}
\end{equation*}
for all $\lambda,\mu\in \Lambda$. So it suffices to prove that
\begin{equation}\label{spec1}
\begin{split}
B\bigl(E(\gamma_{\lambda};z),
\widetilde{T}_i\widetilde{E}(x_{\mu};\widetilde{z})\bigr)
&=\tilde{t}_{i}E(\gamma_{\lambda};x_{\mu}^{-1})\\
&+\tilde{t}_{i}^{-1}
\widetilde{v}_{a_i}(\gamma_{\lambda}^{-1})
\bigl(E(\gamma_{s_i.\lambda};x_{\mu}^{-1})-
E(\gamma_{\lambda};x_{\mu}^{-1})\bigr)
\end{split}
\end{equation}
for all $\lambda,\mu\in \Lambda$ and all $i=0,\ldots,n$.

Formula \eqref{spec1} is easy when $\lambda$ is stabilized by $s_i$
since then we have
\[B\bigl(E(\gamma_{\lambda};z),
\widetilde{T}_i\widetilde{E}(x_{\mu};\widetilde{z})\bigr)=
B\bigl((T_i(E(\gamma_{\lambda};.))(z),\widetilde{E}(x_{\mu};\widetilde{z})\bigr)=
\tilde{t}_iE(\gamma_{\lambda};x_{\mu}^{-1})
\]
where the last equality follows from 
(the proof of) proposition \ref{symmnon} and \eqref{fgexp}.
So we assume for the remainder of the proof
that $s_i.\lambda\not=\lambda$.
We can use now \eqref{extra}
to commute $\widetilde{T}_i$ and
$\widetilde{E}(x_{\mu};\widetilde{z})$
in the left-hand side of \eqref{spec1}. Combined with
lemma \ref{easy} we then derive that
\begin{equation*}
\begin{split}
B\bigl(E(\gamma_{\lambda};z),
\widetilde{T}_i\widetilde{E}(x_{\mu};\widetilde{z})\bigr)&=
\tilde{t}_i\bigl(s_i\widetilde{E}(x_{\mu};\cdot)\bigr)(\gamma_{\lambda}^{-1})\\
&+ 
\psi_i(\gamma_{\lambda}^{-1})\bigl(\widetilde{E}(x_{\mu};\gamma_{\lambda}^{-1})-
(s_i\widetilde{E}(x_{\mu};\cdot))(\gamma_{\lambda}^{-1})\bigr),
\end{split}
\end{equation*}
where $\psi_i\in {\mathbb{C}}(x)$ is given by
\[\psi_i(x)=\frac{(\tilde{t}_{a_i}-\tilde{t}_{a_i}^{-1})+
(\tilde{t}_{a_i/2}-\tilde{t}_{a_i/2}^{-1})x^{a_i/2}}{1-x^{a_i}}=
\tilde{t}_i-\tilde{t}_i^{-1}\widetilde{v}_{a_i}(x).
\]
Since $s_i.\lambda\not=\lambda$, we can apply lemma
\ref{compatibleactions} together with the duality \eqref{duality}
of the non-symmetric Koornwinder polynomials to obtain
the desired formula \eqref{spec1}.
\end{proof}


\section{(Bi-)orthogonality relations and quadratic norms}

{}From now on we assume that $0<q,t<1$ and that the parameters $a,b,c,d$
(see \eqref{parametrization}) have moduli less than one.
We define $\Delta(\cdot)=\Delta(\cdot;\mathbf{t};q)$ and
$\Delta_{+}(\cdot)=\Delta_+(\cdot;\mathbf{t};q)$ by
\begin{equation}
\Delta(x;\mathbf{t};q)=\prod_{\beta\in R^+}\frac{1}{v_\beta(x;\mathbf{t};q)}
\end{equation}
and
\begin{equation}\label{deltaplus}
\begin{split}
\Delta_+(x;\mathbf{t};q)=&\prod_{
\stackrel{\beta\in R}{\beta(0)\geq 0}}
\frac{1}{v_\beta(x;\mathbf{t};q)}\\
=&\prod_{1\leq i<j\leq n}\frac{\bigl(x_ix_j, x_ix_j^{-1},
x_i^{-1}x_j, x_i^{-1}x_j^{-1};q\bigr)_{\infty}}
{\bigl(t^2x_ix_j, t^2x_ix_j^{-1}, t^2x_i^{-1}x_j,
t^2x_i^{-1}x_j^{-1};q\bigr)_{\infty}}\\
&.\prod_{i=1}^n\frac{\bigl(x_i^{2}, x_i^{-2};q\bigr)_{\infty}}
{\bigl(ax_i, ax_i^{-1},
bx_i, bx_i^{-1}, cx_i, cx_i^{-1}, dx_i, dx_i^{-1};q\bigr)_{\infty}},
\end{split}
\end{equation}
where 
$\bigl(y_1,\ldots,y_m;q\bigr)_{\infty}=\prod_{j=1}^m\bigl(y_j;q\bigr)_{\infty}$
with $\bigl(y;q\bigr)_{\infty}=\prod_{j=0}^{\infty}(1-yq^j)$ the $q$-shifted
factorial. The second equality in \eqref{deltaplus} follows from the
${\mathcal{W}}$-orbit structure of the reduced affine
root system $R$ (cf. \eqref{orbit}), together with \eqref{parametrization}.
Observe that $\Delta_+(\cdot)$ is $W$-invariant, where $W$ acts by
permutations and inversions of the coordinates
$x=(x_1,\ldots,x_n)$. Furthermore,
\begin{equation}\label{linkmeasure}
\Delta(x;\mathbf{t};q)= {\mathcal{C}}(x;\mathbf{t};q)\Delta_+(x;\mathbf{t};q)
\end{equation}
with ${\mathcal{C}}(x)={\mathcal{C}}(x;\mathbf{t};q)$ given by
\begin{equation}\label{CC}
{\mathcal{C}}(x;\mathbf{t};q)=\prod_{\alpha\in \Sigma^-}v_{\alpha}(x;\mathbf{t};q).
\end{equation}
We define now 
bilinear forms $\langle .,. \rangle=\langle .,. \rangle_{\mathbf{t},q}$
and $\langle .,. \rangle_+=\langle .,. \rangle_{+,\mathbf{t},q}$
on ${\mathcal{A}}$ by
\begin{equation}
\begin{split}
\langle f,g\rangle&=\frac{1}{(2\pi i)^n}
\iint_{x\in {\mathbb{T}}^n}f(x)(\sigma g)(x)\Delta(x)\frac{dx}{x},\\
\langle f,g\rangle_+&=\frac{1}{(2\pi i)^n}
\iint_{x\in {\mathbb{T}}^n}f(x)(\sigma g)(x)\Delta_+(x)\frac{dx}{x}
\end{split}
\end{equation}
where $\frac{dx}{x}=\frac{dx_1}{x_1}\cdots\frac{dx_n}{x_n}$
and ${\mathbb{T}}\subset {\mathbb{C}}$
is the (positively oriented) unit circle. Recall here that $\sigma$
is the longest Weyl group element in $W$. 
Observe that the bilinear forms $\langle .,. \rangle$
and $\langle .,. \rangle_+$ are non-degenerate in both factors.

The bilinear form $\langle .,. \rangle_+$
coincides with Koornwinder's \cite{Ko} pairing for the
symmetric Koornwinder polynomials, see also \cite{S} for an extension
to more general parameter values. If on the other hand
the parameter values are such that 
$\Delta(\cdot)\in {\mathcal{A}}$, then $\langle f,g \rangle$ equals
the constant term of the Laurent polynomial $f(x)(\sigma g)(x)\Delta(x)$
and $\langle .,. \rangle$ then coincides with Sahi's \cite{Sa2} pairing
for the non-symmetric Koornwinder polynomials.

\begin{lem} We have
\begin{equation}
\label{constant}
\sum_{w\in W}w{\mathcal{C}}(\cdot;\mathbf{t};q)=K_{\mathbf{t},q}
\end{equation}
in ${\mathbb{C}}(x)$ for some constant 
$K=K_{\mathbf{t},q}\in {\mathbb{C}}$. In particular,
\begin{equation}\label{symmred}
\langle f,g\rangle=\frac{K}{|W|}
\langle f,g\rangle_{+},\qquad \forall f,g\in
{\mathcal{A}}^W
\end{equation}
where $|W|=2^nn!$ is the cardinality of the finite Weyl group $W$.
\end{lem}
\begin{proof}
The first statement follows from \cite[(2.8 n.r)]{M2} with the
indeterminates in \cite[(2.8 n.r)]{M2}
specialized to $u_{\alpha}^{1/2}=-t_{\alpha}t_{\alpha/2}$ for 
$\alpha\in \Sigma_l^+$,
$u_{\alpha}=t_{\alpha}^2$ for $\alpha\in \Sigma_m^+$ and
$u_{\alpha}=t_{\alpha}^{-2}$ for $\alpha\in \frac{1}{2}\Sigma_l^+$.
The identity \eqref{symmred} follows then from \eqref{linkmeasure}
and the invariance of the measure $({\mathbb{T}}^n,\frac{dx}{x})$
under the action of $W$.
\end{proof}
A product form for the constant $K$
can be obtained by specializing the left hand side of
\eqref{constant} at $x_0^{-1}$, see \cite[(2.4 n.r)]{M2}.
In fact, we have the following
more general result.
\begin{lem}\label{okpar}
Let $\lambda\in\Lambda^+$, and write $W_{\lambda}$ \textup{(}respectively
$W^{\lambda}$\textup{)} for the stabilizer sub-group of $\lambda$ in $W$
\textup{(}respectively the minimal coset representatives of
$W/W_{\lambda}$\textup{)}. Then
$K=
\sum_{w\in W^{\lambda}}{\mathcal{C}}(x_{w\lambda}^{-1})$. 
In particular,
$K={\mathcal{C}}(x_0^{-1})$.
\end{lem}
\begin{proof}
By the definition \eqref{constant} of $K$ we have
\begin{equation}
\label{sum}
K=\sum_{w\in W^{\lambda}, u\in
W_{\lambda}}(u^{-1}w^{-1}{\mathcal{C}})(x_{\lambda}^{-1}).
\end{equation}
We consider a term $(u^{-1}w^{-1}{\mathcal{C}})(x_{\lambda}^{-1})$ in this
sum with $u\not=1$. Then there exists a simple root $a_i$ ($i\in
\{1,\ldots,n\}$) which is orthogonal to $\lambda$, and which is mapped
to a negative root $\alpha$ by $wu$.
Now remark \ref{casestable} implies that
the factor $v_{u^{-1}w^{-1}\alpha}(x_{\lambda}^{-1})=v_{a_i}(x_{\lambda}^{-1})$
of $(u^{-1}w^{-1}{\mathcal{C}})(x_{\lambda}^{-1})$ is zero.
Hence the contribution in the sum \eqref{sum} is zero unless $u=1$.
The lemma follows now from lemma \ref{compatibleactions}.
\end{proof}

\begin{prop}\label{adjointprop}
For $X\in {\mathcal{H}}$ we have
\[\langle X(f),g\rangle=\langle f, X^{\ddagger}(g)\rangle,\qquad
f,g\in {\mathcal{A}},
\]
where $\ddagger: {\mathcal{H}}\rightarrow {\mathcal{H}}^\prime$ is
the anti-algebra isomorphism defined in lemma \ref{adj}.
\end{prop}
\begin{proof}
The proposition is obviously correct for $X=z^{\lambda}$ ($\lambda\in
\Lambda$), so it suffices to prove it for $X=T_i$ ($i=0,\ldots,n$).
Let $f,g\in {\mathcal{A}}$.
It follows by direct computations that
\begin{equation}\label{adj1}
\bigl(T_if\bigr)(x)\bigl(\sigma
g\bigr)(x)-f(x)\bigl(\sigma\bigl((T_i^\prime)^{-1}g\bigr)\bigr)(x)
=t_i^{-1}h_i(x)v_{a_i}(x)
\end{equation}
for $i=0,\ldots,n$, with
\[h_i(x)=\bigl(s_if\bigr)(x)\bigl(\sigma g\bigr)(x)-
f(x)\bigl(s_i(\sigma g)\bigr)(x)
\]
and with the action of $s_i$ as defined in \eqref{acA}.
Now observe that $h_i$ is $s_i$-alternating, i.e. $s_ih_i=-h_i$ for
$i=0,\ldots,n$. On the other hand,
\begin{equation}\label{fi}
v_{a_i}(x)\Delta(x)=\prod_{\beta\in R^+\setminus \{a_i\}}\frac{1}{v_\beta(x)}
\end{equation}
is invariant under the action of $s_i$ for $i=0,\ldots,n$, where
the action of $s_i$ is extended from ${\mathcal{A}}$ to
(suitably nice) functions $f$ in the $n$ variables $x=(x_1,\ldots,x_n)$
via the formulas \eqref{acA}.
This 
is an immediate consequence of the well-known fact
that the roots $R^+\setminus \{a_i\}$ are permuted by the simple reflection
$s_i$. Hence $\langle T_if,g\rangle-\langle f,(T_i^\prime)^{-1}g\rangle$
can be rewritten as an integral over 
$\bigl({\mathbb{T}}^n,\frac{dx}{x}\bigr)$ with
$s_i$-alternating integrand for all $i\in \{0,\ldots,n\}$.

Now $\langle T_if,g\rangle-\langle f,(T_i^\prime)^{-1}g\rangle=0$ for
$i=1,\ldots,n$ follows from the fact that the measure $\bigl({\mathbb{T}}^n,
\frac{dx}{x}\bigr)$ is $W$-invariant.
The case $i=0$ is more subtle. The
behaviour of the measure $\bigl({\mathbb{T}}^n, \frac{dx}{x}\bigr)$
under the action of $s_0$ is given by
\[ \iint_{x\in {\mathbb{T}}^n}(s_0h)(x)\frac{dx}{x}=
\int_{y_1\in q{\mathbb{T}}}\iint_{y\in {\mathbb{T}}^{n-1}}
h(y_1,y)\frac{dy_1}{y_1}\frac{dy}{y},
\]
which now implies that
\begin{equation}\label{adj2}
\begin{split}
\langle T_0f,g\rangle&-\langle f, (T_0^\prime)^{-1}g\rangle=\\
=&\frac{1}{2(2\pi i)^n}\int_{y_1\in {\mathbb{T}}-q{\mathbb{T}}}\iint_{y\in
{\mathbb{T}}^{n-1}}t_0^{-1}h_0(y_1,y)v_{a_0}(y_1,y)\Delta(y_1,y)\frac{dy_1}{y_1}
\frac{dy}{y}.
\end{split}
\end{equation}
For fixed $y\in {\mathbb{T}}^{n-1}$, the integrand
in the right-hand side of \eqref{adj2} depends analytically on $y_1\in
\{v\in {\mathbb{C}} \, | \, q\leq |v|\leq 1\}$.
Indeed, by a direct computation using \eqref{parametrization} and
the second expression of $\Delta_+(x)$ in \eqref{deltaplus}, we
see that the $y_1$-dependent factor of $v_{a_0}(y_1,y)\Delta(y_1,y)$
is given by
\begin{equation*}
\begin{split}
&\frac{\bigl(y_1^2, q^2y_1^{-2};q\bigr)_{\infty}}
{\bigl(ay_1,by_1,cy_1,dy_1,qay_1^{-1},qby_1^{-1},qcy_1^{-1},
qdy_1^{-1};q\bigr)_{\infty}}\\
&\qquad\qquad\qquad
\qquad\quad.\prod_{j=2}^n\frac{\bigl(y_1y_j, y_1y_j^{-1}, qy_1^{-1}y_j,
qy_1^{-1}y_j^{-1};q\bigr)_{\infty}}
{\bigl(t^2y_1y_j, t^2y_1y_j^{-1}, qt^2y_1^{-1}y_j,
qt^2y_1^{-1}y_j^{-1};q\bigr)_{\infty}},
\end{split}
\end{equation*}
which has the desired analytic behaviour due to the conditions on
the parameters $q$ and $\mathbf{t}$. Thus by Cauchy's theorem we conclude that
$\langle T_0f,g\rangle-\langle f, (T_0^\prime)^{-1}g\rangle=0$.
This completes the proof of the proposition.
\end{proof}

\begin{rem}
An algebraic proof of proposition \ref{adjointprop} was given
by Sahi \cite[thm. 16]{Sa2} for those (discrete) values of
$\mathbf{t}$ such that $\Delta(\cdot;\mathbf{t};q)\in {\mathcal{A}}$.
\end{rem}

We write $E^\prime({\gamma_{\lambda}^{-1}};\cdot)$ for the renormalized
Koornwinder polynomial of degree $\lambda\in \Lambda$
with respect to inverse parameters $(\mathbf{t}^{-1}, q^{-1})$.
Since $(Y^{\lambda})^\ddagger=(Y^{\prime\,\lambda})^{-1}$
for $\lambda\in \Lambda$, we obtain the following extension of
\cite[cor. 17]{Sa2} from theorem \ref{nonsymmetric}
and proposition \ref{adjointprop}.
\begin{cor}[Bi-orthogonality relations]\label{bior}
For $\lambda,\mu\in \Lambda$ with $\lambda\not=\mu$ we have
$\langle E({\gamma_{\lambda}};\cdot),
E^\prime({\gamma_{\mu}^{-1}};\cdot)\rangle=0$.
\end{cor}

Recall that
$E^+({\gamma_\lambda};\cdot)$ ($\lambda\in \Lambda^+$) can be characterized as
the unique solution of the eigenvalue equation
\[Lf=
\bigl(m_{\epsilon_1}(\gamma_{\lambda})-m_{\epsilon_1}(\gamma_0)\bigr)f, \qquad
f\in {\mathcal{A}}^W
\]
which takes the value one at $x_0^{\pm 1}$, where $L$ \eqref{L}
is Koornwinder's second order $q$-difference operator. Since $L$
and the eigenvalue $m_{\epsilon_1}(\gamma_{\lambda})-m_{\epsilon_1}(\gamma_0)$
are invariant under replacement of the parameters $(\mathbf{t},q)$ by
their inverses $(\mathbf{t}^{-1},q^{-1})$, we derive that
\begin{equation}\label{inversesym}
E^+({\gamma_{\lambda}};x;\mathbf{t};q)=E^+({\gamma_{\lambda}^{-1}};x;
\mathbf{t}^{-1};q^{-1}),\qquad \lambda\in \Lambda^+.
\end{equation}
Combined with \eqref{symmred}, corollary \ref{bior} and theorem
\ref{antisymmetric}, we re-obtain Koornwinder's \cite{Ko} orthogonality
relations for the symmetric Koornwinder polynomials:
\begin{cor}[Orthogonality relations]\label{or}
For $\lambda,\mu\in \Lambda^+$ with $\lambda\not=\mu$, we have
$\langle E^+({\gamma_{\lambda}};\cdot), E^+({\gamma_{\mu}};\cdot)\rangle_+=0$.
\end{cor}
We write 
$\hbox{Spec}(Y^\prime)=\{\gamma_{\lambda}^{-1}\, | \, \lambda\in \Lambda\}$
for the spectrum of the $Y^\prime$-operators, and $F=F_{\mathbf{t},q}$
for the linear space of functions 
$g: \hbox{Spec}(Y^\prime)\rightarrow {\mathbb{C}}$
with finite support. We define a ${\mathcal{W}}$-module structure on $F$
by
\[ (wg)(\gamma_{\lambda}^{-1})=g(\gamma_{w^{-1}.\lambda}^{-1}),\qquad
g\in F,\,\, w\in {\mathcal{W}},\,\, \lambda\in \Lambda.
\]

\begin{Def}
We call the linear map 
${\mathcal{F}}={\mathcal{F}}_{\mathbf{t},q}: {\mathcal{A}}
\rightarrow F$ defined by
\begin{equation}
{\mathcal{F}}(f)(\gamma)=\langle f,E^\prime({\gamma};\cdot)\rangle,\qquad
f\in {\mathcal{A}},\,\, \gamma\in \hbox{Spec}(Y^\prime)
\end{equation}
the non-symmetric Koornwinder transform.
\end{Def}
Observe that ${\mathcal{F}}$ is injective since $\langle .,. \rangle$
is non-degenerate, and that ${\mathcal{F}}$ is surjective by
corollary \ref{bior}. 

In the next proposition we present an action of
the double affine Hecke algebra $\widetilde{\mathcal{H}}$ 
on $F$ in terms of spectral difference-reflection operators, and we relate
it to the action of ${\mathcal{H}}$ on ${\mathcal{A}}$ via the 
non-symmetric Koornwinder transform ${\mathcal{F}}$ and
the duality isomorphism $\Phi$.

\begin{prop}\label{F}
The applications
\begin{equation}\label{actdual}
\begin{split}
(\widetilde{T}_ig)(\gamma)&=
\tilde{t}_{i}g(\gamma)+\tilde{t}_{i}^{-1}
\widetilde{v}_{a_i}(\gamma)((s_ig)\bigl(\gamma)-g(\gamma)\bigr),
\qquad i\in \{0,\ldots,n\},\\
\bigl(f(\widetilde{z})g\bigr)(\gamma)&=f(\gamma)g(\gamma),\qquad
f\in {\mathcal{A}}
\end{split}
\end{equation}
where $g\in F$ and $\gamma\in \hbox{Spec}(Y^\prime)$,
uniquely extend to an action of $\widetilde{{\mathcal{H}}}$ on $F$.
Furthermore,
\begin{equation}\label{intertF}
 {\mathcal{F}}(X(f))=\Phi(X){\mathcal{F}}(f),\qquad X\in
{\mathcal{H}},\,\, f\in {\mathcal{A}},
\end{equation}
where $\Phi$ is the duality isomorphism.
\end{prop}
\begin{proof}
The intertwining property \eqref{intertF} can be proved by
checking it for the algebraic generators $U_n$, $T_i$ and $Y_i$
($i=1,\ldots,n$) of ${\mathcal{H}}$ using proposition
\ref{adjointprop} and proposition \ref{transfer}.
The fact that \eqref{actdual} defines an action of $\widetilde{{\mathcal{H}}}$
on $F$ follows then from \eqref{intertF} since $\Phi$ is an
algebra isomorphism and ${\mathcal{F}}$ is bijective.
\end{proof}
Next we determine the inverse of the non-symmetric Koornwinder transform
${\mathcal{F}}$. We let 
${\mathcal{G}}={\mathcal{G}}_{\mathbf{t},q}: F\rightarrow
{\mathcal{A}}$ be the linear map defined by
\begin{equation}
({\mathcal{G}}g)(x)=\sum_{\lambda\in \Lambda}g(\gamma_{\lambda}^{-1})
E({\gamma_{\lambda}};x;\mathbf{t};q)w(\gamma_{\lambda}^{-1};
\tilde{\mathbf{t}};q),\qquad g\in F,
\end{equation}
where the discrete weights $\widetilde{w}(\gamma_{\lambda}^{-1})=
w(\gamma_{\lambda}^{-1};\tilde{\mathbf{t}};q)$ ($\lambda\in\Lambda$) 
are defined as follows:
\begin{equation}
\label{wwplus}
w(\gamma_{\lambda}^{-1};\tilde{\mathbf{t}};q)=
{\mathcal{C}}(\gamma_{\lambda}^{-1};\tilde{\mathbf{t}};q)
w_+(\gamma_{\lambda^+}^{-1};
\tilde{\mathbf{t}};q),
\end{equation}
with $\widetilde{w}_+(\gamma_{\mu}^{-1})=
w_+(\gamma_{\mu}^{-1};\tilde{\mathbf{t}};q)$
for $\mu\in \Lambda^+$ given by the multiple residue
\begin{equation}\label{wplus}
w_+(\gamma_{\mu}^{-1};
\tilde{\mathbf{t}};q)=
\underset{x_1=\gamma_{\mu}^{-\epsilon_1}}{\hbox{Res}}\left(
\underset{x_2=\gamma_{\mu}^{-\epsilon_2}}{\hbox{Res}}\left(
\cdots \underset{x_n=\gamma_{\mu}^{-\epsilon_n}}{\hbox{Res}}
\left(\frac{\Delta_+(x;\tilde{\mathbf{t}};q)}{x_1\cdots x_n}\right)\cdots\right)\right).
\end{equation}
Using the second expression in \eqref{deltaplus} together with \eqref{positive},
it is easily verified that the discrete weights 
$\widetilde{w}(\gamma)$ and $\widetilde{w}_+(\gamma)$
($\gamma\in \hbox{Spec}(Y^\prime)$) are well defined and non-zero
for generic parameters $\mathbf{t}$ and $q$ (in fact,
all residues in \eqref{wplus} are taken at simple poles).

\begin{prop}\label{G}
We have
\[ {\mathcal{G}}(\widetilde{X}g)=\Phi^{-1}(\widetilde{X}){\mathcal{G}}(g),
\qquad \widetilde{X}\in \widetilde{\mathcal{H}},\,\, g\in F.
\]
\end{prop}
\begin{proof}
The proof for $\widetilde{X}=\widetilde{z}^{\lambda}$
with $\lambda\in \Lambda$ is immediate.
Hence it suffices to check the intertwining property for $\widetilde{X}=
\widetilde{T}_i$ ($i=0,\ldots,n$). Let $g\in F$.
By proposition \ref{transfer} we have for $i=0,\ldots,n$,
\[
{\mathcal{G}}\bigl(\widetilde{T}_ig\bigr)-\Phi^{-1}(\widetilde{T}_i)
\bigl({\mathcal{G}}g\bigr)=\tilde{t}_i^{-1}\sum_{\lambda\in\Lambda}
h_i(\gamma_{\lambda};\cdot)\widetilde{v}_{a_i}(\gamma_{\lambda}^{-1})
\widetilde{w}(\gamma_{\lambda}^{-1})
\]
with $h_i(\gamma_{\lambda};\cdot)\in {\mathcal{A}}$ given by
\[
h_i(\gamma_{\lambda};\cdot)=
g(\gamma_{s_i.\lambda}^{-1})E(\gamma_{\lambda};\cdot)
-g(\gamma_{\lambda}^{-1})E(\gamma_{s_i.\lambda};\cdot).
\]
Since $h_i(\gamma_{s_i.\lambda};\cdot)=-h_i(\gamma_{\lambda};\cdot)$ for
$i=0,\ldots,n$
and $\lambda\in\Lambda$, it thus suffices to prove that
\begin{equation}\label{lastbit}
\widetilde{v}_{a_i}(\gamma_{\lambda}^{-1})
\widetilde{w}(\gamma_{\lambda}^{-1})=
\widetilde{w}_+(\gamma_{\lambda^+}^{-1})
\prod_{\alpha\in \Sigma^-\cup
\{a_i\}}\widetilde{v}_{\alpha}(\gamma_{\lambda}^{-1})
\end{equation}
is invariant under replacement of $\lambda\in \Lambda$ by $s_i.\lambda$
for all $i\in \{0,\ldots,n\}$ and all $\lambda\in \Lambda$.
For $i\in\{1,\ldots,n\}$ this is immediate by lemma
\ref{compatibleactions}.

As usual, the proof for the affine part of the statement 
(the case $i=0$) is more subtle.
We begin by rewriting 
$\widetilde{w}(\gamma_{\lambda}^{-1})$
as a (kind of) multiple residue of 
$\widetilde{\Delta}(x)$
at $x=\gamma_{\lambda}^{-1}$, where $\widetilde{\Delta}(x)=
\Delta(x;\tilde{\mathbf{t}};q)$.
This can be done
using the $W$-invariance of the weight function
$\Delta_+(\cdot;\tilde{\mathbf{t}};q)$. The result is as follows.

Let $u_{\lambda}\in S_n$ be the component in $S_n$ of
the minimal coset representative $w_{\lambda}\in W$ with respect to the
semi-direct product structure $W=S_n\ltimes (\pm 1)^n$, and let
$n_{\lambda}=\#\{ i\in \{1,\ldots,n\} \, | \, \lambda_i<0\}$.
Then we have
\begin{equation}\label{Resnieuw}
\widetilde{w}(\gamma_{\lambda}^{-1})=
\underset{x=\gamma_{\lambda}^{-1}}{\hbox{\bf{Res}}}\left(
\frac{\widetilde{\Delta}(x)}{x_1\cdots
x_n}\right)
\end{equation}
for all $\lambda\in\Lambda$, 
where the multiple residue at $x=\gamma_{\lambda}^{-1}$
is defined by
\[
\underset{x=\gamma_{\lambda}^{-1}}{\hbox{\bf{Res}}}\bigl(\,\cdot\,\bigr)=
(-1)^{n_{\lambda}}
\underset{x_{u_{\lambda}(1)}=\gamma_{\lambda}{}^{-\epsilon_{u_{\lambda}(1)}}}
{\hbox{Res}}
\Bigl(\underset{x_{u_{\lambda}(2)}=
\gamma_{\lambda}{}^{-\epsilon_{u_{\lambda}(2)}}}
{\hbox{Res}}\Bigl(\cdots \underset{x_{u_{\lambda}(n)}=
\gamma_{\lambda}{}^{-\epsilon_{u_{\lambda}(n)}}}{\hbox{Res}}\,\Bigl(\cdot\Bigr)
\cdots\Bigr)\Bigr).
\]
In particular, we obtain
\begin{equation}\label{altexpression}
\widetilde{v}_{a_0}(\gamma_{\lambda}^{-1})
\widetilde{w}(\gamma_{\lambda}^{-1})=
\underset{x=\gamma_{\lambda}^{-1}}{\hbox{\bf{Res}}}\left(
\frac{\widetilde{v}_{a_0}(x)\widetilde{\Delta}(x)}{x_1\cdots
x_n}\right)
\end{equation}
for all $\lambda\in\Lambda$.
Now we consider \eqref{altexpression} with $\lambda$ replaced by
$s_0.\lambda$. We first consider the changes in the multiple
residue. By the proof of \cite[thm. 5.3]{Sa}, we have
$w_{s_0.\lambda}=s_{\epsilon_1}w_{\lambda}$ for all
$\lambda\in\Lambda$, i.e.
$n_{s_0.\lambda}=n_{\lambda}\pm 1$ and $u_{s_0.\lambda}=u_{\lambda}$. Secondly,
\begin{equation*}
(\gamma_{s_0.\lambda}^{-1})^{\epsilon_i}=
\begin{cases}
q\gamma_{\lambda}^{\epsilon_1}\qquad &\hbox{ if } i=1,\\
\gamma_{\lambda}^{-\epsilon_i} &\hbox{ if } i=2,\ldots,n
\end{cases}
\end{equation*}
by lemma \ref{compatibleactions}. So if we replace
the residue at $x_1=\gamma_{\lambda}^{-\epsilon_1}$ by the residue
at $x_1=q\gamma_{\lambda}^{\epsilon_1}$ in the definition of the
multiple residue at $x=\gamma_{\lambda}^{-1}$, then we obtain minus the
multiple residue at $x=\gamma_{s_0.\lambda}^{-1}$.
On the other hand, we know by the proof of proposition
\ref{adjointprop} that $v_{a_0}(x)\Delta(x)$  is invariant under
the action of $s_0$. So the invariance of \eqref{altexpression}
under replacement of $\lambda$
by $s_0.\lambda$ follows from the simple observation that
\[ \underset{y=y_0}{\hbox{Res}}\left(\frac{h(y)}{y}\right)=
-\underset{y=qy_0^{-1}}{\hbox{Res}}\left(\frac{h(y)}{y}\right)
\]
when $h(y)$ is a function depending on a single variable $y$,
having a simple pole at
$y=y_0$, and satisfying the invariance condition
$h(qy^{-1})=h(y)$.
\end{proof}

\begin{thm}\label{mainr}
We have ${\mathcal{G}}\circ {\mathcal{F}}=k\, \hbox{Id}_{\mathcal{A}}$
and ${\mathcal{F}}\circ {\mathcal{G}}=k\,\hbox{Id}_{F}$
with $k=k_{\mathbf{t},q}=w(\gamma_0^{-1};\tilde{\mathbf{t}};q)\,
\langle 1,1\rangle_{\mathbf{t},q}$.
In particular, we have for $\lambda,\mu\in \Lambda$,
\begin{equation}\label{norm1}
\frac{\langle E({\gamma_{\lambda}};\cdot),
E^\prime(\gamma_{\mu}^{-1};\cdot)\rangle}
{\langle 1,1\rangle}=
\delta_{\lambda,\mu}
\frac{\widetilde{w}(\gamma_0^{-1})}{\widetilde{w}(\gamma_{\lambda}^{-1})}
\end{equation}
where $\delta_{\lambda,\mu}$ is the Kronecker delta.
\end{thm}
\begin{proof}
By proposition \ref{F} and proposition \ref{G} we have
\begin{equation}\label{ok1}
{\mathcal{G}}({\mathcal{F}}(f))={\mathcal{G}}({\mathcal{F}}(f(z)1))=
f(z){\mathcal{G}}({\mathcal{F}}(1)),\qquad \forall f\in
{\mathcal{A}}.
\end{equation}
Furthermore, it follows from corollary
\ref{bior} that
\begin{equation}
\label{ok2}
{\mathcal{G}}\bigl({\mathcal{F}}(E({\gamma};\cdot))\bigr)=
\langle E({\gamma};\cdot),
E^\prime({\gamma^{-1}};\cdot)\rangle\,
\widetilde{w}(\gamma^{-1})\,E({\gamma};\cdot)
\end{equation}
for $\gamma\in \hbox{Spec}(Y)$. Formula \eqref{ok2} reduces
to ${\mathcal{G}}({\mathcal{F}}(1))=k\,1$ when
$\gamma=\gamma_0$, with the constant $k$ as given in the statement
of the theorem. Combined with \eqref{ok1} it follows
that ${\mathcal{G}}\circ {\mathcal{F}}=
k\,\hbox{Id}_{\mathcal{A}}$. Since ${\mathcal{F}}$ is bijective,
we then also have ${\mathcal{F}}\circ {\mathcal{G}}=k\,\hbox{Id}_F$.

It remains to prove \eqref{norm1}.
By corollary \ref{bior}, we only have to prove \eqref{norm1}
when $\mu=\lambda$.
We fix $\gamma=\gamma_{\lambda}\in \hbox{Spec}(Y)$,  $\lambda\in \Lambda$.
Since ${\mathcal{G}}\circ
{\mathcal{F}}=k\,\hbox{Id}_{\mathcal{A}}$,
it follows that ${\mathcal{G}}({\mathcal{F}}(E(\gamma;\cdot)))=
k\,E(\gamma;\cdot)$. Comparing this outcome with
the right-hand side of \eqref{ok2}, we obtain
\[\langle E({\gamma};\cdot),
E^\prime({\gamma^{-1}};\cdot)\rangle\,
\widetilde{w}(\gamma^{-1})=k=\langle
1,1\rangle\, \widetilde{w}(\gamma_0^{-1})\]
which yields the desired result.

\end{proof}

\begin{cor}\label{normplus}
For all $\lambda,\mu\in \Lambda^+$, we have
\[
\frac{\langle E^+({\gamma_{\lambda}};\cdot),
E^+({\gamma_{\mu}};\cdot)\rangle_{+}}
{\langle 1,1\rangle_+}=
\delta_{\lambda,\mu}
\frac{\widetilde{w}_+(\gamma_0^{-1})}{\widetilde{w}_+(\gamma_{\lambda}^{-1})}.
\]
\end{cor}
\begin{proof}
First of all, observe that $E^+({\gamma_{\lambda^+}};\cdot)=
C_+E({\gamma_{\lambda}};\cdot)$
for all $\lambda\in \Lambda$ and that
$(C_+)^\ddagger=C_+^\prime$, where
$C_+\in {\mathcal{H}}$ (respectively $C_+^\prime\in
{\mathcal{H}}^\prime$, $\widetilde{C}_+ \in \widetilde{\mathcal{H}}$) 
is the idempotent corresponding to the trivial representation
of the underlying finite Hecke algebra of type $C_n$. Combined with
proposition \ref{adjointprop}, \eqref{symmred} and
\eqref{inversesym} we can rewrite ${\mathcal{F}}_+={\mathcal{F}}|_{{\mathcal{A}}^W}$
as
\[ \bigl({\mathcal{F}}_+f\bigr)(\gamma_{\lambda}^{-1})=
\frac{K}{|W|}\,\langle f,
E^+({\gamma_{\lambda^+}};\cdot)\rangle_+,\qquad f\in
{\mathcal{A}}^W,\,\, \lambda\in\Lambda.
\]
In particular, the symmetric Koornwinder transform ${\mathcal{F}}_+$ maps into
\[
F^W=\{g\in F \, | \, wg=g\,\,\, \forall w\in {\mathcal{W}}\}
=\{ g\in F \, | \, \widetilde{C}_+g=g \}.
\]
Similarly, since $\Phi(C_+)=\widetilde{C}_+$, we derive from proposition \ref{G},
\eqref{wwplus} and lemma \ref{okpar} that
${\mathcal{G}}_+={\mathcal{G}}|_{F^W}$ can be rewritten as
\[\bigl({\mathcal{G}}_+g\bigr)(x)=
\widetilde{K}\sum_{\lambda\in
\Lambda^+}g(\gamma_{\lambda}^{-1})E^+({\gamma_{\lambda}};x)
\widetilde{w}_+(\gamma_{\lambda}^{-1}),
\]
where $\widetilde{K}=K_{\tilde{\mathbf{t}},q}$.
Using these alternative descriptions for ${\mathcal{F}}_+$ and ${\mathcal{G}}_+$
together with the orthogonality relations for the symmetric
Koornwinder polynomials (see corollary \ref{or}), we obtain
\[ {\mathcal{G}}_+({\mathcal{F}}_+(E^+(\gamma_{\lambda};\cdot)))=
\frac{K\widetilde{K}}{|W|}\,
\langle E^+(\gamma_{\lambda};\cdot),
E^+(\gamma_{\lambda};\cdot)\rangle_{+}\,
\widetilde{w}_+(\gamma_{\lambda}^{-1})\,
E^+(\gamma_{\lambda};\cdot)
\]
for all $\lambda\in\Lambda^+$. 
On the other hand, by theorem
\ref{mainr}, 
we have ${\mathcal{G}}_+({\mathcal{F}}_+(E^+(\gamma_{\lambda};\cdot)))=
k\,E^+(\gamma_{\lambda};\cdot)$ for $\lambda\in\Lambda^+$ with
\[k=\langle
1,1\rangle\,\widetilde{w}(\gamma_0^{-1})=
\frac{K\widetilde{K}}{|W|}\,
\langle 1,1\rangle_{+}\,\widetilde{w}_+(\gamma_0^{-1}).
\]
Comparing the two different outcomes for
${\mathcal{G}}_+({\mathcal{F}}_+(E^+(\gamma_{\lambda};\cdot)))$,
we obtain the desired result in case $\mu=\lambda$. The
off-diagonal case is covered by the orthogonality relations
for the symmetric Koornwinder polynomials, see corollary \ref{or}.
\end{proof}

\begin{rem}\label{explicweight}
The discrete weights 
$\widetilde{w}_+(\gamma_{\lambda}^{-1})=w_+(\gamma_{\lambda}^{-1};
\tilde{\mathbf{t}};q)$
($\lambda\in\Lambda^+$) also appear as weights in a
partly discrete orthogonality measure for 
the symmetric Koornwinder polynomials, see \cite{S}. 
In particular, using the expression \cite[prop. 4.1]{S}
for the discrete weights, it can be shown that the ratio
$\langle E^+(\gamma_{\lambda};\cdot), E^+(\gamma_{\lambda};\cdot)\rangle_+/
\langle 1,1\rangle_+=\widetilde{w}_+(\gamma_0^{-1})/
\widetilde{w}_+(\gamma_{\lambda}^{-1})$ is equal to 
\begin{equation*}
\begin{split}
&\prod_{i=1}^n\left\{\frac{\bigl(q^{-1}abcdt^{4(n-i)};q\bigr)_{2\lambda_i}
\bigl(c^2t^{4(n-i)}\bigr)^{\lambda_i}}
{\bigl(abcdt^{4(n-i)};q\bigr)_{2\lambda_i}}\right.\\
&\qquad\qquad\qquad\qquad
\qquad\qquad\left..\frac{\bigl(qt^{2(n-i)}, abt^{2(n-i)}, adt^{2(n-i)},
bdt^{2(n-i)};q\bigr)_{\lambda_i}}
{\bigl(act^{2(n-i)}, bct^{2(n-i)}, cdt^{2(n-i)},
q^{-1}abcdt^{2(n-i)};q\bigr)_{\lambda_i}}\right\}\\
&.\prod_{1\leq i<j\leq n}\frac{\bigl(abcdt^{2(2n-i-j-1)},
q^{-1}abcdt^{2(2n-i-j)};q\bigr)_{\lambda_i+\lambda_j}}
{\bigl(abcdt^{2(2n-i-j)},
q^{-1}abcdt^{2(2n-i-j+1)};q\bigr)_{\lambda_i+\lambda_j}}
\frac{\bigl(qt^{2(j-i-1)},
t^{2(j-i)};q\bigr)_{\lambda_i-\lambda_j}}
{\bigl(qt^{2(j-i)}, t^{2(j-i+1)};q\bigr)_{\lambda_i-\lambda_j}}
\end{split}
\end{equation*}
for all $\lambda\in\Lambda^+$, where $\bigl(y_1,\ldots,y_m;q\bigr)_k=
\prod_{j=1}^m\bigl(y_j;q\bigr)_k$ and 
$\bigl(y;q\bigr)_k=\prod_{j=0}^{k-1}(1-yq^j)$
for $k\in {\mathbb{Z}}_+$. 
\end{rem}


\section{Evaluation formulas}

In this section we give an explicit expression
for the value of the non-symmetric Koornwinder polynomial 
$P_{\lambda}(\cdot)=P_{\lambda}(\cdot;\mathbf{t};q)$ at $x_0^{-1}$. 
We start with two preliminary lemmas.

\begin{lem}\label{lem1}
For $\lambda\in \Lambda$, we have
\[\frac{\langle P_{\lambda}, E^\prime({\gamma_{\lambda}^{-1}};\cdot)\rangle}
{\langle P_{\lambda^+},
E^\prime({\gamma_{\lambda^+}^{-1}};\cdot)\rangle}
=\tilde{t}_{w_{\lambda}}^2\prod_{\alpha\in \Sigma^+\cap
w_{\lambda}^{-1}\Sigma^-}
\frac{1}{\widetilde{v}_{\alpha}(\gamma_{\lambda^+}^{-1})}.
\]
\end{lem}
\begin{proof}
We prove the lemma using the non-affine intertwiners $S_i$
($i=1,\ldots,n$). For the moment, we fix $\lambda\in \Lambda$
and $i\in\{1,\ldots,n\}$ such that $\mu=s_i\lambda\not=\lambda$.
We first give some additional properties of $S_i$ which we will need
for the proof.

The intertwiner $S_i$ is self-adjoint, i.e. $(S_i)^\ddagger=S_i^\prime$, where
$S_i^\prime=[T_i^\prime,(Y^\prime)^{a_i}]\in
{\mathcal{H}}^\prime$.
This can be checked most easily in the image of the duality
anti-isomorphism $\Psi_{\mathbf{t}^{-1},q^{-1}}$. It follows from
proposition \ref{transfer}{\bf (ii)} that 
$S_i^\prime\,E^\prime({\gamma_{\lambda}^{-1}};\cdot)=
L_{i,\lambda}E^\prime({\gamma_{\mu}^{-1}};\cdot)$ with 
$L_{i,\lambda}=
\bigl(\gamma_{\lambda}^{-a_i}-\gamma_{\lambda}^{a_i}\bigr)
\tilde{t}_{i}^{-1}\widetilde{v}_{a_i}(\gamma_{\lambda}^{-1})$.
On the other hand, $S_iP_{\lambda}=K_{i,\lambda}P_{\mu}$ with
$K_{i,\lambda}=\bigl(\gamma_{\lambda}^{a_i}-\gamma_{\lambda}^{-a_i}\bigr)
\eta_i(\gamma_\lambda)$ by corollary \ref{normintertwiner}. Combined with
proposition \ref{adjointprop} we obtain
\[
\langle P_{\mu}, E^\prime({\gamma_{\mu}^{-1}};\cdot)\rangle
=\frac{1}{L_{i,\lambda}}\langle S_iP_{\mu},
E^\prime({\gamma_{\lambda}^{-1}};\cdot)\rangle
=
\frac{K_{i,\mu}}{L_{i,\lambda}}\langle P_{\lambda},
E^\prime({\gamma_{\lambda}^{-1}};\cdot)\rangle.
\]
In particular, if $\langle \lambda, a_i\rangle<0$, then
\[
\langle P_{\lambda}, E^\prime({\gamma_{\lambda}^{-1}};\cdot)\rangle=
\frac{\tilde{t}_{i}^2}{\widetilde{v}_{a_i}(\gamma_{\mu}^{-1})}
\langle P_{\mu}, E^\prime({\gamma_{\mu}^{-1}};\cdot)\rangle.
\]
The ratio
$\langle P_{\lambda}, E^\prime({\gamma_{\lambda}^{-1}};\cdot)\rangle/
\langle P_{\lambda^+}, E^\prime({\gamma_{\lambda^+}^{-1}};\cdot)\rangle$
can now be evaluated inductively using similar techniques as in the proof of
theorem \ref{antisymmetric}. This gives the desired result.
\end{proof}

Let $\delta_{\gamma_{\lambda}^{-1}}\in F$ for $\lambda\in\Lambda$
be the function which is equal to one at $\gamma_{\lambda}^{-1}$
and zero otherwise.

\begin{lem}\label{lem2}
For $\lambda\in \Lambda^+$ we have
\[\bigl(\Phi(z^{\lambda})\delta_{\gamma_0^{-1}}\bigr)
(\gamma_{\lambda}^{-1})=\tilde{t}_{\tau(-\lambda)}^{-1}
\prod_{\beta\in R^+\cap \tau(\lambda)R^-}
\widetilde{v}_\beta(\gamma_{\lambda}^{-1}).
\]
\end{lem}
\begin{proof}
Let $\lambda\in \Lambda^+$. It follows from \eqref{length} that
$\tau(-\lambda)$ is the unique element of minimal length
in ${\mathcal{W}}$ which maps $\lambda$ to $0\in \Lambda$ under the dot-action.
In particular, any element $w\in {\mathcal{W}}$ which is smaller than
$\tau(-\lambda)$ with respect to the Bruhat order, maps $\lambda$
to a non-zero element in $\Lambda$.

Since $\lambda\in \Lambda^+$, we have $Y^{\lambda}=T_{\tau(\lambda)}$, hence
\[ \Phi(z^{\lambda})=(\epsilon(z^{\lambda}))^{\widetilde{\dagger}^\prime}=
(\widetilde{T}^\prime_{\tau(\lambda)})^{\widetilde{\dagger}^\prime}=
\widetilde{T}^{-1}_{\tau(-\lambda)}.
\]
Let now $\tau(-\lambda)=s_{i_1}s_{i_2}\cdots s_{i_r}$ be
a reduced expression of $\tau(-\lambda)$ in ${\mathcal{W}}$, then
we obtain from proposition \ref{F} and from the previous
paragraph that
\[
\bigl(\Phi(z^{\lambda})\delta_{\gamma_0^{-1}}\bigr)(\gamma_{\lambda}^{-1})
=\prod_{m=1}^r\tilde{t}_{a_{i_m}}^{-1}
\widetilde{v}_{a_{i_m}}\bigl(\gamma_{(s_{i_{m+1}}\cdots
s_{i_{r-1}}s_{i_r}).\lambda}^{-1}\bigr),
\]
where $(s_{i_{m+1}}\cdots s_{i_{r-1}}s_{i_r}).\lambda$ should be read
as $\lambda$ when $m=r$. The lemma follows now easily from
lemma \ref{compatibleactions} and remark \ref{compatible}.
\end{proof}

\begin{thm}\label{evaluationthm}
Let $\lambda\in \Lambda$, then
\[ P_{\lambda}(x_0^{-1})=
\frac{\widetilde{w}(\gamma_{\lambda}^{-1})}
{\widetilde{w}(\gamma_0^{-1})}\,\tilde{t}_{w_{\lambda}}^2
\tilde{t}_{\tau(-\lambda^+)}^{-1}\prod_{\beta}
\widetilde{v}_\beta(\gamma_{\lambda^+}^{-1}),
\]
where the product is taken over $\beta\in (R^+\cap \tau(\lambda^+)R^-)\setminus
(\Sigma^+\cap w_{\lambda}^{-1}\Sigma^-)$.
\end{thm}
\begin{proof}
First of all, observe that
\begin{equation}\label{severalforms}
\begin{split}
P_{\lambda}(x_0^{-1})\langle E({\gamma_{\lambda}};\cdot),
E^\prime({\gamma_{\lambda}^{-1}};\cdot)\rangle&=
\langle P_{\lambda}, E^\prime({\gamma_{\lambda}^{-1}};\cdot)\rangle,\\
\langle P_{\lambda}, E^\prime({\gamma_{\lambda}^{-1}};\cdot)\rangle
&=\langle x^{\lambda}, E^\prime({\gamma_{\lambda}^{-1}};\cdot)\rangle=
\langle 1,1\rangle\bigl(\Phi(z^{\lambda})\delta_{\gamma_0^{-1}}\bigr)
(\gamma_{\lambda}^{-1})
\end{split}
\end{equation}
for all $\lambda\in \Lambda$. Indeed, the first formula is trivial,
while the second formula is a direct consequence of 
corollary \ref{bior}, the definition of the non-symmetric Koornwinder
transform ${\mathcal{F}}$ and its intertwining properties given in 
proposition \ref{F}.

We use now successively the first formula of \eqref{severalforms},
then lemma \ref{lem1} and finally the second formula of \eqref{severalforms}
to arrive at
\[ P_{\lambda}(x_0^{-1})=
\frac{\tilde{t}_{w_{\lambda}}^2
\langle 1,1\rangle}{\langle E(\gamma_{\lambda};\cdot),
E^\prime(\gamma_{\lambda}^{-1};\cdot)\rangle}
\bigl(\Phi(z^{\lambda^+})\delta_{\gamma_0^{-1}}\bigr)(\gamma_{\lambda^+}^{-1})
\prod_{\alpha\in\Sigma^+\cap
w_{\lambda}^{-1}\Sigma^-}\frac{1}
{\widetilde{v}_{\alpha}(\gamma_{\lambda^+}^{-1})}.
\]
The theorem follows now from theorem \ref{mainr} and lemma \ref{lem2},
combined with the inclusion $\Sigma^+\cap w_{\lambda}^{-1}\Sigma^-\subset
R^+\cap\tau(\lambda^+)R^-$. This inclusion is a direct consequence
of the inequality $\langle \lambda^+,\alpha\rangle> 0$
for all $\alpha\in \Sigma^+\cap w_{\lambda}^{-1}\Sigma^-$, see \eqref{char}.
\end{proof}
\begin{cor}\label{evalplus}
\begin{equation}\label{evalform}
P_{\lambda}^+(x_0^{\pm 1})=
\frac{\widetilde{w}_+\bigl(\gamma_{\lambda}^{-1})}
{\widetilde{w}_+\bigl(\gamma_0^{-1})}\,
\tilde{t}_{\tau(-\lambda)}^{-1}\prod_{\beta\in
R^+\cap\tau(\lambda)R^-}
\widetilde{v}_\beta\bigl(\gamma_{\lambda}^{-1}),
\qquad \lambda\in \Lambda^+.
\end{equation}
\end{cor}
\begin{proof}
This follows from the evaluation of $P_{\mu}(x_0^{-1})$ ($\mu\in
\Lambda$) (see theorem \ref{evaluationthm}), theorem
\ref{antisymmetric}{\bf (iii)}, \eqref{wwplus} and lemma
\ref{okpar}.
\end{proof}

\begin{rem}\label{evalrem}
It is straightforward to explicitly write down
the roots $R^+\cap \tau(\lambda)R^-$ for $\lambda\in\Lambda^+$.
Together with \eqref{positive}, \eqref{v}, remark \ref{explicweight}
and the ${\mathcal{W}}$-orbit structure of $R$, one can
reformulate \eqref{evalform} now as
\begin{equation*}
\begin{split}
P_{\lambda}^+(x_0)=&\prod_{i=1}^n
\frac{\bigl(act^{2(n-i)}, bct^{2(n-i)}, cdt^{2(n-i)},
q^{-1}abcdt^{2(n-i)};q\bigr)_{\lambda_i}}
{\bigl(q^{-1}abcdt^{4(n-i)};
q\bigr)_{2\lambda_i}\bigl(ct^{2(n-i)}\bigr)^{\lambda_i}}\\
&.\prod_{1\leq i<j\leq
n}\frac{\bigl(q^{-1}abcdt^{2(2n-i-j+1)};q\bigr)_{\lambda_i+\lambda_j}
\bigl(t^{2(j-i+1)};q\bigr)_{\lambda_i-\lambda_j}}
{\bigl(q^{-1}abcdt^{2(2n-i-j)};q\bigr)_{\lambda_i+\lambda_j}
\bigl(t^{2(j-i)};q\bigr)_{\lambda_i-\lambda_j}}
\end{split}
\end{equation*}
for $\lambda\in\Lambda^+$.
\end{rem}

\begin{rem}
Van Diejen \cite[thm. 5.1]{vD2} proved the
evaluation formula \eqref{evalform} for a five parameter sub-family of the 
symmetric Koornwinder polynomials. This result was indirectly extended
to the complete six parameter family of symmetric Koornwinder polynomials by
Sahi's \cite{Sa} duality results.
\end{rem}


\section{The generalized Weyl character formula and the constant term}

In this section we discuss several results involving
the anti-symmetric Koornwinder polynomials $P_{\lambda}^-(\cdot)=
P_{\lambda}^-(\cdot;\mathbf{t};q)$ ($\lambda\in\Lambda^{++}$) in an essential
way. The most crucial result is the analogue of the Weyl character
formula for the Koornwinder polynomials.

In order to state the generalized Weyl character formula, we 
need to introduce some notations first.
Let ${\mathbf{q}}=\{q_\beta\}_{\beta\in S}$ be the multiplicity
function satisfying $q_{a_0}=q_{a_0^\vee}=q_{a_n^\vee}=1$,
$q_{a_i}=q^{1/2}$ ($i\in \{1,\ldots,n-1\}$) and $q_{a_n}=q$.
For two multiplicity functions $\mathbf{t}$ and $\mathbf{t}^\prime$,
we write $\mathbf{t}\mathbf{t}^\prime$ for the multiplicity
function which takes the value $t_{\beta}t_{\beta}^\prime$ at $\beta\in S$.
Then the generalized Weyl character formula 
is given by
\begin{equation}\label{gWcf}
P_{\lambda+\kappa}^-\bigl(x;\mathbf{t};q)=
\chi(x;\mathbf{t};q)P_{\lambda}^+(x;\mathbf{q}\mathbf{t};q),\qquad \lambda
\in\Lambda^+
\end{equation}
with $\kappa$ given by \eqref{kappa} and
with $\chi(\cdot;\mathbf{t};q)\in {\mathcal{A}}$ given by
\[
\chi(x;\mathbf{t};q)=x^{\kappa}\prod_{\alpha\in\Sigma^-}
(1-x^{\alpha})v_{\alpha}(x;\mathbf{t}^{-1};q^{-1}).
\]
The proof of \eqref{gWcf}
is similar to the proof of the generalized Weyl character formula for
Macdonald polynomials, see e.g. \cite[7.3]{M3}.

The generalized Weyl character formula \eqref{gWcf}, together with
the results of section 8, readily implies the norm relations
\begin{equation}\label{normrelations}
\begin{split}
\frac{\langle P_{\lambda}^+(\cdot;\mathbf{t};q), P_{\lambda}^+(\cdot;
\mathbf{t};q)\rangle_{\mathbf{t},q}}
{\langle P_{\lambda}^-(\cdot;\mathbf{t};q), P_{\lambda}^-(\cdot;
\mathbf{t}^{-1};q^{-1})\rangle_{\mathbf{t},q}}&=
t_{\sigma}^2\frac{\langle P_{\lambda}^+(\cdot;\mathbf{t};q),
P_{\lambda}^+(\cdot;\mathbf{t};q)\rangle_{+,\mathbf{t},q}}
{\langle P_{\lambda-\kappa}^+(\cdot;\mathbf{q}\mathbf{t};q),
P_{\lambda-\kappa}^+(\cdot;\mathbf{q}\mathbf{t};q)\rangle_{+,
\mathbf{q}\mathbf{t},q}}\\
&=\prod_{\alpha\in\Sigma^+}
\frac{v_{\alpha}(\gamma_{\lambda}^{-1};\tilde{\mathbf{t}};q)}
{v_{\alpha}(\gamma_{\lambda};\tilde{\mathbf{t}};q)}
\end{split}
\end{equation}
for $\lambda\in\Lambda^{++}$. The norm relations \eqref{normrelations}
give an explicit description of the diagonal terms 
$\langle P_{\lambda}^-(\cdot;\mathbf{t};q), 
P_{\lambda}^-(\cdot;\mathbf{t}^{-1};q^{-1})\rangle_{\mathbf{t},q}$
($\lambda\in\Lambda^{++}$)
corresponding to the bi-orthogonality relations
\begin{equation}\label{bioranti} 
\langle P_{\lambda}^-(\cdot;\mathbf{t};q),
P_{\mu}^-(\cdot;\mathbf{t}^{-1};q^{-1})\rangle_{\mathbf{t},q}=0\qquad
\lambda,\mu\in\Lambda^{++}: \,\lambda\not=\mu
\end{equation}
in terms of the quadratic
norms of the symmetric Koornwinder polynomials. 
Furthermore, the norm relations \eqref{normrelations} enables one 
to obtain a new proof for Gus\-taf\-son's \cite{G} 
evaluation of the constant term $\langle 1,1 \rangle_+$.
In fact, by \eqref{normrelations}, one can express $\langle 1,1\rangle_{+,
\mathbf{q}^m\mathbf{v},q}$ in terms of 
\begin{equation}\label{step}
\langle P_{m\kappa}^+(\cdot;\mathbf{v};q),P_{m\kappa}^+(\cdot;\mathbf{v};q)
\rangle_{+,\mathbf{v},q}
\end{equation}
for all positive integers $m$. Let now $\mathbf{v}$ be 
a multiplicity function with $v_\beta=1$ for all $\beta\in S$ of length two,
then the corresponding orthogonality measure
reduces to the (coordinate-wise) product measure 
of the one-variable Askey-Wilson polynomials. In particular, \eqref{step}
can be expressed in terms of the quadratic norms of the 
Askey-Wilson polynomials, which were evaluated in  
\cite{AW} (see \cite{NS} for an affine
Hecke algebraic approach). This yields an evaluation of
\eqref{step}, and hence of $\langle 1,1\rangle_{+,\mathbf{q}^m\mathbf{v},q}$.
By analytic continuation, we arrive at Gustafson's \cite{G} result that
the constant term $|W|^{-1}\langle 1,1\rangle_{+,\mathbf{t},q}$ is equal
to
\[\prod_{j=1}^n
\frac{\bigl(t^2,t^{2(2n-j-1)}abcd;q\bigr)_{\infty}}
{\bigl(q,t^{2(n-j+1)}, t^{2(n-j)}ab, t^{2(n-j)}ac, t^{2(n-j)}ad,
t^{2(n-j)}bc, t^{2(n-j)}bd, t^{2(n-j)}cd;q\bigr)_{\infty}}.
\]
\begin{rem}
In view of theorem \ref{mainr}, corollary \ref{normplus}, 
theorem \ref{evaluationthm}, corollary \ref{evalplus} and
Gustafson's constant term evaluation, we have arrived now
at the stage that the quadratic norms (respectively diagonal terms)
of the (non-)symmetric Koornwinder polynomials are completely
explicit. In particular, the explicit evaluation of $\langle P_{\lambda}^+, 
P_{\lambda}^+\rangle_+$ ($\lambda\in\Lambda^+$) which we thus obtain, 
can be seen to
coincide with van Diejen's \cite[thm. 5.2]{vD2} explicit expression for
$\langle P_{\lambda}^+, P_{\lambda}^+\rangle_+$.
\end{rem}


\bibliographystyle{amsplain}

\end{document}